\documentclass[MSNbibl,number,citesort,dvips]{arxbj}
\usepackage{upgreek}
\usepackage{graphicx}

\aid{0}
\volume{19}
\issue{2}
\pubyear{2013}
\firstpage{462}
\lastpage{491}
\doi{10.3150/12-BEJ413} 

\makeatletter
\newcommand{\eqref}[1]{(\ref{#1})}
\newcommand{\B}{\boldsymbol}
\newcommand{\E}{\mathrm{E}}
\newcommand{\Var}{\operatorname{Var}}
\newcommand{\Cov}{\operatorname{Cov}}

\newcommand{\R}{\mathbb{R}}
\newcommand{\N}{\mathcal{N}}
\newcommand{\I}{\mathcal{I}}
\newcommand{\J}{\mathcal{J}}
\newcommand{\K}{\mathcal{K}}
\newcommand{\Q}{\mathcal{Q}}
\newcommand{\Rr}{\mathcal{R}}

\newremark{Ex}{Example}[section]
\newtheorem{Theorem}{Theorem}
\newproclaim{Cond}{Condition}

\makeatother

\begin{document}
\begin{frontmatter}

\title{Parameter estimation for pair-copula constructions}
\runtitle{Parameter estimation for pair-copula constructions}

\begin{aug}
\author{\fnms{Ingrid} \snm{Hob{\AE}k Haff}\corref{}\ead[label=e1]{ingrid@nr.no}}
\runauthor{I. Hob{\ae}k Haff} 
\address{Statistics for Innovation, Norwegian Computing Center, PB 114 Blindern, NO-0373 Oslo,
\mbox{Norway}. \printead{e1}}
\end{aug}

\received{\smonth{12} \syear{2010}}
\revised{\smonth{9} \syear{2011}}

%
\begin{abstract}
We explore various estimators for the parameters
of a pair-copula construction (PCC), among those the stepwise
semiparametric (SSP) estimator, designed for this dependence
structure. We present its asymptotic properties, as well as the
estimation algorithm for the two most common types of PCCs.
Compared to the considered alternatives, that is, maximum likelihood,
inference functions for margins and semiparametric estimation,
SSP is in general asymptotically less efficient. As we show in
a few examples, this loss of efficiency may however be rather low.
Furthermore, SSP is semiparametrically efficient for the Gaussian
copula. More importantly, it is computationally tractable even
in high dimensions, as opposed to its competitors. In any case, SSP
may provide start values, required by the other estimators. It is
also well suited for selecting the pair-copulae of a PCC for a
given data set.
\end{abstract}

%

\end{frontmatter}

\section{Introduction}\label{sec:introduction}
The last decades' technological revolution have considerably
increased the relevance of multivariate modelling. Copulae are
now regularly used within fields such as finance, survival
analysis and actuarial sciences. Although the list of parametric
bivariate copulae is long and varied, the choice is rather limited
in higher dimensions (Genest \textit{et~al.}~\cite{VinesEd}).
Accordingly, a number of
hierarchical, copula-based structures have been proposed, among
those the pair-copula construction (PCC) of Joe~\cite{joe96}, further
studied and considered by Bedford and Cooke
\cite{bedfordcooke01,bedfordcooke02}, Kurowicka and Cooke
\cite{kurowickacooke06} and Aas \textit{et~al.}~\cite{Vines}.

\begin{figure}

\includegraphics{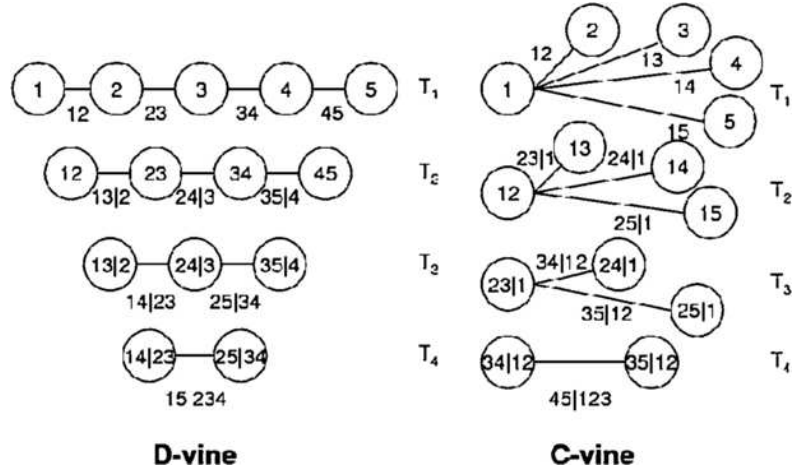}

\caption{Five-dimensional D-vine (to the left) and C-vine (to the
right).}\label{fig:CDvine}
\end{figure}

A PCC is a treelike construction, built from pair-copulae with
conditional distributions as their two arguments (see Figure
\ref{fig:CDvine}). The number of conditioning variables is zero
at the ground, and increases by one for each level, to ensure
coherence of the construction. Despite its simple structure, the PCC
is highly flexible and covers a wide range of complex dependencies (Joe
\textit{et~al.}~\cite{joe09}, Hob{\ae}k Haff \textit{et~al.}~\cite{PCCcover}).
After Aas \textit{et~al.}~\cite{Vines} set it in an inferential
context, it has made several appearances in the literature (Fischer \textit{et~al.}
\cite{fischer07}, Chollete \textit{et~al.}~\cite{valdesogo08}, Heinen and Valdesogo
\cite{valdesogo082}, Schirmacher and Schirmacher~\cite{schirmacher08}, Czado and Min~\cite{czadomin08},
Kolbj{\o}rnsen and Stien~\cite{kolbjstien08}, Czado \textit{et~al.}~\cite{czadomin09,czadoschepmin10}),
exhibiting its adequacy for various applications.

Regardless of its recent popularity, estimation of PCC
parameters has so far been addressed mostly in an applied
setting. The aim of this work is to explore the properties
of alternative estimators. As the PCC is a member of the
multivariate copula family, one may exploit the large
collection of estimators proposed for that model class, such as
moments type procedures, based on, for instance, the matrix
of pairwise Kendall's tau  coefficients (Clayton~\cite{clayton78}, Oakes~\cite{oakes82},
Genest~\cite{genest87}, Genest and Rivest~\cite{genest93}). Such methods may be well-suited for particular
copula families. We are however interested in more general
procedures, allowing for broader model classes. Moreover, we
wish to exploit the specific structure of the PCC.

More specifically, the number of parameters of a PCC grows
quickly with the dimension, even if all pair-copulae
constituting the structure are from one-parameter families.
In medium to high dimension, the existing copula estimators
may simply become too demanding computationally, and will at
least require good start values in the optimisation procedure.
Furthermore, due to the PCC's tree structure, selection of
appropriate pair-copulae for a given data set must be done
level by level. Procedures that estimate all parameters
simultaneously are therefore unfit for this task.

In all, we contemplate four estimators. The first is the
classical maximum likelihood (ML), followed by the inference
functions for margins (IFM) and semiparametric estimators, that
have been developed specifically for multivariate copulae. These three
estimators are treated in Section~\ref{sec:multicopest}, and are
included mostly for comparison. Section~\ref{sec:pccest} is
devoted to the fourth one, the stepwise semiparametric estimator
(SSP). Unlike the others, it is designed for the PCC structure.
Although it has been suggested and used earlier (Aas \textit{et~al.}~\cite{Vines}), it
has never been formally presented,\vadjust{\goodbreak} nor have its asymptotic properties
been explored. In Section~\ref{sec:examples}, we compare the four
estimators in a few examples. Finally, Section~\ref{sec:conclusion}
presents some concluding remarks.

The setting is as follows. Consider the observations
$\mathbf{x}_{1},\ldots,\mathbf{x}_{n}$ of $n$ independent
$d$-variate stochastic vectors $\mathbf{X}_{1},\ldots,\mathbf{X}_{n}$,
originating from the same pair-copula construction. Assume
further that the joint distribution is absolutely continuous,
with strictly increasing margins. The corresponding copula is
then unique (Sklar~\cite{sklar}). Letting $\B{\alpha}$ and $\B{\theta}$
denote the parameters of the margins and copula, respectively,
the joint probability density function (p.d.f.) may then be
expressed as (McNeil \textit{et~al.}~\cite{mcneil06}, page 197)\vspace*{-1.5pt}
%
\begin{equation}
f_{1\ldots d}(x_{1},\ldots,x_{d};\B{\alpha},\B{\theta})  =
c_{1\ldots d}(F_{1}(x_{1};\B{\alpha}_{1}),\ldots,F_{d}(x_{d};\B{\alpha
}_{d});\B{\theta})\prod_{l=1}^{d}f_{l}(x_{l};\B{\alpha}_{l}).
\label{eqn:genpdf}
\end{equation}
Here, $F_{l}$ and $f_{l}, l=1,\ldots,d$, are the marginal
cumulative distribution functions (c.d.f.s) and p.d.f.s, respectively,
and $c_{1\ldots d}$ is the corresponding copula density. Since this is a
PCC, $c_{1\ldots d}$ is, in turn, a product of pair-copulae.

\begin{table}
\tablewidth=\textwidth
\tabcolsep=0pt
  \caption{Notation overview}\label{tab:notation}
\begin{tabular*}{\textwidth}{@{\extracolsep{\fill}}llll@{}}
  \hline
  Symbol & $j = 0$ & $j = 1$ & $j = 2,\ldots,d-1$ \\
  \hline
  $v_{ij}^{*}$ & -- & $\varnothing$ & $\{i+1,\ldots,i+j-1\}$ \\
  $w_{ij}^{*}$ & $i$ & $\{i,i+j\}$ & $\{i,v_{ij},i+j\}$ \\
  $\B{\theta}_{j}^{*}$ & -- & $\{\B{\theta}_{i,i+1}\dvt i = 1,\ldots,d-1\}$ & $\{\B{\theta}_{s,s+t|v_{st}}\dvt|v_{st}| = j-1\}$ \\
  $\B{\theta}_{i\rightarrow i+j}^{*}$ & -- & $\B{\theta}_{i,i+1}$ & $\{\B{\theta}_{s,s+t|v_{st}}\dvt (s,s+t) \in w_{ij}\}$ \\
  $z_{ij}^{**}$ & $\varnothing$ & $\{i+1\}$ & $\{1,\ldots,j-1,j+i\}$ \\
  $\B{\theta}_{j\_i}^{**}$ & $\varnothing$ & $\varnothing$ & $\{\B{\theta}_{s,s+t|z_{s-1,0}}\dvt (s,s+t) \in z_{ji}\}$ \\
  \hline
\end{tabular*}
  \tabnotetext[]{}{$^{*}$: $i=1,\ldots,d-j$, $^{**}$:
  $i=0,\ldots,d-j$.}\vspace*{-1.5pt}
\end{table}

Define the index sets $v_{ij}=\{i+1,\ldots,i+j-1\}$,
$w_{ij}=\{i,v_{ij},i+j\}$, for $1\leq i \leq d-j, 1\leq j \leq d-1$,
with $v_{i1} = \varnothing$, and $w_{i0}=i$. Thus, for a
vector $\mathbf{a}=(a_{1},\ldots,a_{d})$, we write
$\mathbf{a}_{v_{ij}}=(a_{i+1},\ldots,a_{i+j-1})$ and
$\mathbf{a}_{w_{ij}}=(a_{i},\ldots,a_{i+j})$. Further, for an
index $k$ and a set of indices $v$, with $k \notin v$, let
$F_{k|v}$ be the conditional c.d.f. of $X_{k}$ given
$\mathbf{X}_{v}=\mathbf{x}_{v}$, and $c_{i,i+j|v_{ij}}$ the
copula density corresponding to the conditional distribution
$F_{i,i+j|v_{ij}}$ of $(X_{i},X_{i+j})$ given
$\mathbf{X}_{v_{ij}}=\mathbf{x}_{v_{ij}}$. Finally, let
$\B{\theta}_{i,i+j|v_{ij}}$ be the parameters of the
copula density $c_{i,i+j|v_{ij}}$, and define
$\B{\theta}_{i\rightarrow i+j} = \{\B{\theta
}_{s,s+t|v_{st}}\dvt(s,s+t) \in w_{ij}\}$,
with $\B{\theta}_{i\rightarrow i} = \varnothing$, and
$\B{\theta}_{j} = \{\B{\theta}_{s,s+t|v_{st}}\dvt |v_{st}| = j-1\}$,
where $|\cdot|$ denotes the cardinality (i.e., $\B{\theta}_{j}$
gathers all parameters from level $j$ of the structure). Table
\ref{tab:notation} gives an overview of the notation. For a
so-called D-vine (Bedford and Cooke~\cite{bedfordcooke01,bedfordcooke02}), the
joint p.d.f. \eqref{eqn:genpdf} can now be written as
(Aas \textit{et~al.}~\cite{Vines})
%
\begin{eqnarray}\label{eqn:dvinepdf}
&&f_{1\ldots d}(x_{1},\ldots,x_{d};\B{\alpha},\B{\theta})\nonumber\\
&&\quad = \prod_{l=1}^{d}f_{l}(x_{l};\B{\alpha}_{l}) \prod_{j=1}^{d-1}\prod_{i=1}^{d-j}c_{i,i+j|v_{ij}}
(F_{i|v_{ij}}(x_{i}|\mathbf{x}_{v_{ij}};
\B{\alpha}_{w_{i,j-1}},\B{\theta}_{i\rightarrow i+j-1}) , \\
&&\quad \hphantom{= \prod_{l=1}^{d}f_{l}(x_{l};\B{\alpha}_{l}) \prod_{j=1}^{d-1}\prod_{i=1}^{d-j}c_{i,i+j|v_{ij}}
(} F_{i+j|v_{ij}}(x_{i+j}|\mathbf{x}_{v_{ij}};\B{\alpha}_{w_{i+1,j-1}},
\B{\theta}_{i+1\rightarrow i+j});\B{\theta}_{i,i+j|v_{ij}} ).\nonumber
\end{eqnarray}
In four dimensions, this becomes
%
\begin{eqnarray}\label{eqn:dvinepdf4D}
f_{1234}(x_{1},x_{2}, x_{3},x_{4};\B{\alpha},\B{\theta})
&=&
f_{1}(x_{1};\B{\alpha}_{1})\cdot f_{2}(x_{2};
\B{\alpha}_{2})\cdot f_{3}(x_{3};\B{\alpha}_{3})\cdot
f_{4}(x_{4};\B{\alpha}_{4})\nonumber\\
&&\hspace*{-0.5pt}{} \cdot c_{12}(F_{1}(x_{1};\B{\alpha}_{1}),F_{2}(x_{2};
\B{\alpha}_{2});\B{\theta}_{12})\nonumber\\
&&\hspace*{-0.5pt}{}\cdot c_{23}(F_{2}(x_{2};\B{\alpha
}_{2}),F_{3}(x_{3};
\B{\alpha}_{3});\B{\theta}_{23})\nonumber \\
& &\hspace*{-0.5pt}{}\cdot c_{34}(F_{3}(x_{3};\B{\alpha}_{3}),F_{4}(x_{4};\B{\alpha
}_{4});\B{\theta}_{34})\\
&&{} \cdot c_{13|2}(F_{1|2}(x_{1}|x_{2};\B{\alpha}_{1},\B{\alpha
}_{2},\B{\theta}_{12}),
F_{3|2}(x_{3}|x_{2};\B{\alpha}_{2},\B{\alpha}_{3},\B{\theta
}_{23});\B{\theta}_ {13|2})\nonumber\\
&&{} \cdot c_{24|3}(F_{2|3}(x_{2}|x_{3};\B{\alpha}_{2},\B{\alpha
}_{3},\B{\theta}_{23}),
F_{4|3}(x_{4}|x_{3};\B{\alpha}_{3},\B{\alpha}_{4},\B{\theta
}_{34});\B{\theta}_{24|3})\nonumber\\
&&{} \cdot c_{14|23}(F_{1|23}(x_{1}|x_{2},x_{3};\B{\alpha}_{1},\B
{\alpha}_{2},\B{\alpha}_{3},
\B{\theta}_{12},\B{\theta}_{23},\B{\theta}_{13|2}),\nonumber\\
&&\hphantom{{} \cdot c_{14|23}(}F_{4|23}(x_{4}|x_{2},x_{3};\B{\alpha}_{2},\B{\alpha}_{3},\B
{\alpha}_{4},
\B{\theta}_{23},\B{\theta}_{34},\B{\theta}_{24|3});\B{\theta}_{14|23}).\nonumber
\end{eqnarray}

For simplicity, we will start by assuming that the distribution
in question is a D-vine, represented to the left in Figure
\ref{fig:CDvine} for $d=5$. Similar results can be obtained
for C-vines (Section~\ref{subsec:c-vines}) and other
regular vines (Bedford and Cooke~\cite{bedfordcooke01,bedfordcooke02}). We also
assume that the PCC is of a simplified form (Hob{\ae}k Haff \textit{et~al.}~\cite{PCCcover}),
that is, that the parameters $\B{\theta}_{i,i+j|v_{ij}}$ of the
copulae $C_{i,i+j|v_{ij}}$, combining conditional distributions,
are not functions of the conditioning variables $\mathbf{x}_{v_{ij}}$.
Without this assumption, inference on these models is not
doable in practice.
%
\section{Multivariate copula estimators}\label{sec:multicopest}
As previously mentioned, a PCC is a multivariate copula. Hence, one
may estimate its parameters with well-known methods, such as maximum
likelihood or the two-step inference functions for margins and
semiparametric estimators.

\subsection{Maximum likelihood (ML) estimator}\label{subsec:ml}
Supposing the model is true, the ML estimator is a natural choice, due
to its asymptotic efficiency and other advantageous characteristics.
According to \eqref{eqn:dvinepdf}, the log-likelihood function of a D-vine
is given by
%
\begin{eqnarray}\label{eqn:loglik}
l(\B{\alpha},\B{\theta};\mathbf{x}) &=& \sum_{k=1}^{n}\log
(f_{1\ldots d}(x_{1k},\ldots,x_{dk};\B{\alpha},\B{\theta}))\nonumber\\
& = &\sum_{k=1}^{n}\sum_{l=1}^{d}\log(f_{l}(x_{lk};\B{\alpha
}_{l}))\nonumber\\
& &{}+\sum_{k=1}^{n}\sum_{j=1}^{d-1}\sum_{i=1}^{d-j}\log
(c_{i,i+j|v_{ij}} (F_{i|v_{ij}}(x_{i}|\mathbf{x}_{v_{ij}};
\B{\alpha}_{w_{i,j-1}},\B{\theta}_{i\rightarrow i+j-1}),
\\
&&\hphantom{{}+\sum_{k=1}^{n}\sum_{j=1}^{d-1}\sum_{i=1}^{d-j}\log
(c_{i,i+j|v_{ij}} (} F_{i+j|v_{ij}}(x_{i+j}|\mathbf{x}_{v_{ij}};
\B{\alpha}_{w_{i+1,j-1}},\B{\theta}_{i+1\rightarrow i+j});\B
{\theta}_{i,i+j|v_{ij}} ) )\nonumber\\
& =& l_{M}(\B{\alpha};\mathbf{x})+l_{C}(\B{\alpha},\B{\theta};\mathbf{x}),\nonumber
\end{eqnarray}
where $\mathbf{x}=(\mathbf{x}_{1},\ldots,\mathbf{x}_{n})$. The ML estimator
$\hat{\B{\theta}}{}^{\mathrm{ML}}$ is obtained by maximising the above log-likelihood
function over all parameters, $\B{\theta}$ and $\B{\alpha}$, simultaneously.
Under the additional assumptions (M1)--(M8) of Lehmann~\cite{asymp} (pages 499--501),
this corresponds to solving the set of estimating equations (one equation
per parameter),
$\frac{1}{n}\sum_{k=1}^{n}\B{\phi}^{\mathrm{ML}} (X_{1k},\ldots,X_{dk};
\hat{\B{\alpha}}{}^{\mathrm{ML}},\hat{\B{\theta}}{}^{\mathrm{ML}} )=\B{0}$,
which
is a vector of functions, with elements
%
\begin{eqnarray}\label{eqn:mlesteqn}
\B{\phi}_{l}^{\mathrm{ML}}(x_{1},\ldots,x_{d};\B{\alpha},\B{\theta}) & =&
\frac{\partial\log(f_{1\ldots d}(x_{1},\ldots,x_{d};
\B{\alpha},\B{\theta}) )}{\partial\B{\alpha}_{l}},\nonumber
\\[-8pt]
\\[-8pt]
\B{\phi}_{d+(j-1) (d-j/2 )+i}^{\mathrm{ML}}(x_{1},\ldots,x_{d};\B
{\alpha},\B{\theta}) & =&
\frac{\partial\log(f_{1\ldots d}(x_{1},\ldots,x_{d};\B{\alpha},\B
{\theta}) )}{\partial\B{\theta}_{i,i+j|v_{ij}}}\nonumber
\end{eqnarray}
for $l = 1,\ldots,d$, $i = 1,\ldots,d-j$, $j = 1,\ldots,d-1$. Define
$\B{\I}$ as the corresponding Fisher information matrix
\begin{eqnarray*}
\B{\I} & =& \E\biggl( \biggl(\frac{\partial\log(f_{1\ldots d}(\mathbf{X};\B{\alpha
},\B{\theta}) )}{\partial(\B{\alpha},
\B{\theta})} \biggr) \biggl(\frac{\partial\log(f_{1\ldots d}(\mathbf{X};\B{\alpha},\B
{\theta}) )}{\partial(\B{\alpha},
\B{\theta})} \biggr)^{T} \biggr)\\
& = &\E\biggl(-\frac{\partial^{2} \log(f_{1\ldots d}(\mathbf{X};\B{\alpha},\B
{\theta}) )}{\partial(\B{\alpha},
\B{\theta})\,\partial(\B{\alpha},\B{\theta})^{T}} \biggr) =
\pmatrix{
\B{\I}_{\alpha} & \B{\I}_{\alpha,\theta}\cr
\B{\I}_{\alpha,\theta}^{T} & \B{\I}_{\theta}
}.
\end{eqnarray*}
In the last expression, it is partitioned according to marginal and
dependence parameters. The corresponding inverse is
%
\begin{eqnarray}\label{eqn:invfisher}
\B{\I}^{-1} & =&\pmatrix{
\B{\I}^{(\alpha)} & \B{\I}^{(\alpha,\theta)}\cr
\bigl(\B{\I}^{(\alpha,\theta)} \bigr)^{T} & \B{\I}^{(\theta)}
},\nonumber\\
\B{\I}^{(\alpha)} &=& (\B{\I}_{\alpha}-\B{\I
}_{\alpha,\theta}\B{\I}_{\theta}^{-1}\B{\I}_{\alpha,\theta
}^{T} )^{-1},\nonumber
\\[-8pt]
\\[-8pt]
\B{\I}^{(\alpha,\theta)}& =& -\B{\I}^{(\alpha)}\B{\I}_{\alpha
,\theta}\B{\I}_{\theta}^{-1},\nonumber\\
\B{\I}^{(\theta)}& =& \B{\I}_{\theta}^{-1}+\B{\I}_{\theta
}^{-1}\B{\I}_{\alpha,\theta}^{T}\B{\I}^{(\alpha)}\B{\I
}_{\alpha,\theta}\B{\I}_{\theta}^{-1}.\nonumber
\end{eqnarray}
It is well-known that under the mentioned conditions, the
estimator $\hat{\B{\theta}}{}^{\mathrm{ML}}$ is consistent for $\B{\theta}$
and asymptotically normal, that is,
\begin{eqnarray*}
\sqrt{n} (\hat{\B{\theta}}{}^{\mathrm{ML}}-\B{\theta} )& \stackrel
{d}{\longrightarrow}& \N(\B{0},\mathbf{V}^{\mathrm{ML}} ),\\
\mathbf{V}^{\mathrm{ML}} &=& \B{\I}^{(\theta)}.
\end{eqnarray*}
In general, ML estimation of PCC parameters will require numerical
optimisation. Even in rather low dimensions, such as four or
five, the number of parameters is high if several of the model
components have more than one parameter. For instance, a
five-dimensional PCC, consisting of Student's \textit{t}-copulae, has 20
parameters, to which one must add the ones of the margins. Finding
the global maximum in such a high-dimensional space is numerically
challenging, even with more elaborate optimisation schemes, such as
the modified Newton--Raphson method with first and second order derivatives.
It will in any case be highly time consuming, so the ML estimator
may not be an option in practice. Therefore, one needs faster and
computationally easier estimation procedures.

Moreover, the above results require the chosen model to be the true
model, that is, the one that produced the data. If the specified
model is close to the truth in the Kullback--Leibler (KL) sense, the
ML estimator may behave very well (Claeskens and Hjort~\cite{claeskens08}). However, it
is in general non-robust to larger KL-divergences from the true model.
\subsection{Two-step estimators}\label{subsec:2step}
The next two estimators are not particularly designed for
pair-copula constructions, but for multivariate copula
models in general. Both consist of two steps, the first
being estimation of the marginal parameters.
\subsubsection{Inference function for margins (IFM) estimator}\label{subsubsec:ifm}
The IFM estimator, introduced by Joe~\cite{joe97,joe05}, addresses
the computational inefficiency of the ML estimator by performing the
estimation in two steps. First, one estimates $\B{\alpha}$ by maximising
the term $l_{M}$ from \eqref{eqn:loglik}. The resulting estimates
$\hat{\B{\alpha}}{}^{\mathrm{IFM}}$ are plugged into the term $l_{C}$ to obtain
$\hat{\B{\theta}}{}^{\mathrm{IFM}}$.

Under conditions (M1)--(M8) (see Section~\ref{subsec:ml}),
this corresponds to solving
\[\frac{1}{n}\sum_{k=1}^{n}\B{\phi}^{\mathrm{IFM}} (X_{1k},\ldots
,X_{dk};\hat{\B{\alpha}}{}^{\mathrm{IFM}},\hat{\B{\theta}}{}^{\mathrm{IFM}} )=\B{0},\]
with elements
\begin{eqnarray}\label{eqn:ifmesteqn}
\B{\phi}_{l}^{\mathrm{IFM}}(x_{l};\B{\alpha}_{l}) & =& \frac{\partial\log
(f_{l}(x_{l};\B{\alpha}_{l}) )}{\partial\B{\alpha}_{l}},\nonumber
\\[-8pt]
\\[-8pt]
\B{\phi}_{d+(j-1) (d-j/2 )+i}^{\mathrm{IFM}}(x_{1},\ldots,x_{d};\B
{\alpha},\B{\theta}) & = &\frac{\partial\log(f_{1\ldots
d}(x_{1},\ldots,x_{d};\B{\alpha},\B{\theta}) )}{\partial\B
{\theta}_{i,i+j|v_{ij}}}\nonumber
\end{eqnarray}
for $l = 1,\ldots,d$, $i = 1,\ldots,d-j$, $j = 1,\ldots,d-1$.
Compared to the
ML equations \eqref{eqn:mlesteqn}, the full log-p.d.f., $\log f_{1\ldots
d}$, is
replaced with the marginal log-p.d.f.s, $\log f_{j}$, for the estimation of
$\B{\alpha}$.

Consider a four-dimensional D-vine \eqref{eqn:dvinepdf4D}, consisting
of Student's \textit{t}-copulae, each having their own correlation and degrees
of freedom parameter, combined with Student's \textit{t}-margins. The parameter
vectors are then $\B{\alpha}=(\nu_{1},\nu_{2},\nu_{3},\nu_{4})$ and
$\B{\theta}=(\rho_{12},\rho_{23},\rho_{34},\rho_{13|2},\rho_{24|3},\rho_{14|23},\nu_{12},\nu_{23},\nu_{34},\nu
_{13|2},\nu_{24|3},\nu_{14|23})$.
IFM estimation of this model starts with a separate estimation
of $\nu_{i}$, $i=1,2,3,4$, margin by margin. The next step is to
optimise
$l_{C} (\hat{\nu}_{1},\ldots,\hat{\nu}_{4},\B{\theta};\mathbf{x} )$
over $\B{\theta}$, $l_{C}$ being the sum of the log-copula densities in
line 3 to 8 of \eqref{eqn:dvinepdf4D}, over all observations.

Define the matrix $\B{\K}_{\alpha}$ with
\begin{eqnarray*}
\B{\K}_{\alpha,i,j} &=&\E
(\B{\phi}_{i}^{\mathrm{IFM}}(X_{i};\B{\alpha}_{i})\B{\phi
}_{j}^{\mathrm{IFM}}(X_{j};\B{\alpha}_{j}) )\\
&=&\E\biggl( \biggl(\frac{\partial\log(f_{i}(X_{i};\B{\alpha}_{i}) )}{\partial
\B{\alpha}_{i}} \biggr) \biggl(\frac{\partial\log(f_{j}(X_{j};\B{\alpha
}_{j}) )}{\partial\B{\alpha}_{j}} \biggr)^{T} \biggr),
\qquad i,j=1,\ldots,d,
\end{eqnarray*}
and the block diagonal matrix $\B{\J}_{\alpha}$ with
\[
\B{\J}_{\alpha,i,i} = \E\biggl(-\frac{\partial}{\partial\B{\alpha
}_{i}^{T}}\B{\phi}_{i}^{\mathrm{IFM}}(X_{i};\B{\alpha}_{i})\biggr ) = \E\biggl(-\frac
{\partial^{2}\log(f_{i}(X_{i};\B{\alpha}_{i}) )}{\partial\B
{\alpha}_{i}\,\partial\B{\alpha}_{i}^{T}} \biggr) = \B{\K}_{\alpha,i,i},
\]
each block corresponding to one of
the margins. If all margins are one-parameter families,
$\B{\K}_{\alpha}$ and $\B{\J}_{\alpha}$ are $d\times d$ matrices.
More generally, their dimension depends on the number of
parameters of each margin.

Joe~\cite{joe05} showed that under the mentioned conditions, the
estimator $\hat{\B{\theta}}{}^{\mathrm{IFM}}$ is consistent for $\B{\theta}$,
as well as asymptotically normal:
%
\begin{eqnarray}\label{eqn:covifm}
 \sqrt{n} (\hat{\B{\theta}}{}^{\mathrm{IFM}}-\B{\theta} ) &\stackrel
{d}{\longrightarrow}& \N(\B{0},\mathbf{V}^{\mathrm{IFM}} ),\nonumber
\\[-8pt]
\\[-8pt]
 \mathbf{V}^{\mathrm{IFM}} &=& \B{\I}_{\theta}^{-1}+\B{\I}_{\theta}^{-1}\B{\I
}_{\alpha,\theta}^{T}\B{\J}_{\alpha}^{-1}\B{\K}_{\alpha}
\B{\J}_{\alpha}^{-1}\B{\I}_{\alpha,\theta}\B{\I}_{\theta}^{-1}.\nonumber
\end{eqnarray}
The above covariance matrix is obtained by replacing $\B{\I}^{(\alpha)}$
in \eqref{eqn:invfisher} with the asymptotic covariance matrix
$\B{\J}_{\alpha}^{-1}\B{\K}_{\alpha}\B{\J}_{\alpha}^{-1}$ of
$\hat{\B{\alpha}}{}^{\mathrm{IFM}}$. This quantifies the loss of asymptotic
efficiency from discarding information the dependence structure might
have on the margins. Several studies, including Joe~\cite{joe05} and
Kim \textit{et~al.}~\cite{silvapulle07}, have demonstrated that unless the dependence
between the variables is extreme, this loss tends to be rather small.
That is also the impression from Examples~\ref{ex:expgumb} and \ref
{ex:student}
(Section~\ref{sec:examples}). Moreover, the IFM method is
computationally faster than the ML estimator, and can at least
be used to set the start values in ML optimisation. Of course,
for high-dimensional $\B{\theta}$, the IFM estimator is still
too slow to be used for PCCs.
\subsubsection{Semiparametric (SP) estimator for copula
parameters}\label{subsubsec:sp}
Just like IFM, the SP estimator is a two-step estimator, treating the
margins separately. It was introduced by Genest \textit{et~al.}~\cite{genest95}, and for the
censored case by Shih and Louis~\cite{shih95}. Later, it was generalised by
Tsukahara~\cite{tsukahara05}. Aas \textit{et~al.}~\cite{Vines} suggest this estimator for PCCs.

As seen from \eqref{eqn:dvinepdf}, the pair-copula arguments at
the ground level of a pair-copula construction (T\tsub{1} in Figure
\ref{fig:CDvine}) are marginal distributions $F_{i}$. From the
second level, they are conditional distributions, whose
conditioning set increases by one with each level. These conditional
distributions may however be written as functions of the margins. Let
$i,j$ be distinct
indices, that is, $i \neq j$, and $v$ a nonempty set of indices, all
from $\{1,\ldots,d\}$, such that $i,j \notin v$. Then, in a
simplified pair-copula construction (Joe~\cite{joe96})
%
\begin{equation}\label{eqn:seqcond}
F_{i|v\cup j}(x_{i}|\mathbf{x}_{v\cup j}) = \frac{\partial
C_{ij|v}(u_{i},u_{j})}{\partial u_{j}} \bigg|_{u_{i}=F_{i|v}(x_{i}|\mathbf
{x}_{v}),u_{j}=F_{j|v}(x_{j}|\mathbf{x}_{v})}.
\end{equation}
Thus, by extracting one of the variables $j$ from the conditioning
set $v\cup j$, one can express $F_{i|v\cup j}$ as a function of
two conditional distributions $F_{i|v}$ and $F_{j|v}$ with one
conditioning variable less. Likewise, $F_{i|v}$ and $F_{j|v}$ may
be written as bivariate functions of conditional distributions with
a conditioning set reduced by one. Proceeding in this way, one finally
obtains recursive functions of the margins.

The type of pair-copula construction determines which conditional
distributions are needed. At level $j \geq2$ of a D-vine, these
are the pairs $(F_{i|v_{ij}}(x_{i}|\mathbf
{x}_{v_{ij}}),F_{i+j|v_{ij}}(x_{i+j}|\mathbf{x}_{v_{ij}})),\break
i=1,\ldots,d-j$. Now, define the functions
%
\begin{eqnarray}\label{eqn:h-func}
h_{i,i+j|v_{ij}}(u_{i},u_{i+j}) & \equiv&\frac{\partial
C_{i,i+j|v_{ij}}(u_{i},u_{i+j})}{\partial u_{i+j}},\nonumber
\\[-8pt]
\\[-8pt]
h_{i+j,i|v_{ij}}(u_{i+j},u_{i}) & \equiv&\frac{\partial
C_{i,i+j|v_{ij}}(u_{i},u_{i+j})}{\partial u_{i}}
\nonumber
\end{eqnarray}
for $i=1,\ldots,d-j, j=1,\ldots,d-1$. Using \eqref{eqn:seqcond},
one obtains
\begin{eqnarray*}
&&F_{i|v_{ij}}(x_{i}|\mathbf{x}_{v_{ij}})\\
&&\quad= h_{i,i+j-1|v_{i,j-1}}(F_{i|v_{i,j-1}}(x_{i}|\mathbf
{x}_{v_{i,j-1}}),F_{i+j-1|v_{i,j-1}}(x_{i+j-1}|\mathbf{x}_{v_{i,j-1}})),\\
&&F_{i+j|v_{ij}}(x_{i+j}|\mathbf{x}_{v_{i,j-1}})\\
 &&\quad=
 h_{i+j,i+1|v_{i+1,j-1}}(F_{i+j|v_{i+1,j-1}}(x_{i+j}|\mathbf
{x}_{v_{i+1,j-1}}),F_{i+1|v_{i+1,j-1}}(x_{i+1}|\mathbf{x}_{v_{i+1,j-1}})),
\end{eqnarray*}
which are functions of conditional distributions constituting
the arguments of the copulae from the previous level, $j-1$. As one
continues this recursion, one achieves, as earlier mentioned, functions
of the margins $F_{i},\ldots,F_{i+j}$. Since these are needed in the
asymptotics, we denote them $g_{i,i+j}^{(1)}$ and $g_{i,i+j}^{(2)}$, and
explicitly define them below. Note however, that for all\vspace*{1pt} practical
purposes, such as in the estimation algorithm (Algorithm~1 in
\ref{sup:SSPalg} of Hob{\ae}k Haff~\cite{ho2012}), one will use the nested $h$-functions from
\eqref{eqn:h-func}. Define
%
\begin{eqnarray}\label{eqn:g-func}
g_{i,i+j}^{(1)}(u_{i},\ldots,u_{i+j-1})
&\equiv&
F_{i|v_{ij}}(F_{i}^{-1}(u_{i})|F_{i+1}^{-1}(u_{i+1}),\ldots
,F_{i+j-1}^{-1}(u_{i+j-1})),\nonumber
\\[-8pt]
\\[-8pt]
g_{i,i+j}^{(2)}(u_{i+1},\ldots,u_{i+j})
&\equiv&
F_{i+j|v_{ij}}(F_{i+j}^{-1}(u_{i+j})|F_{i+1}^{-1}(u_{i+1}),\ldots
,F_{i+j-1}^{-1}(u_{i+j-1}))\nonumber
\end{eqnarray}
for $ i=1,\ldots,d-j, j=1,\ldots,d-1$. Now, one may
rewrite \eqref{eqn:dvinepdf} as:
%
\begin{eqnarray}\label{eqn:dvinepdfg}
&&f_{1\ldots d}(x_{1},\ldots,x_{d};\B{\alpha},\B{\theta})\nonumber\\
 &&\quad=\prod_{l=1}^{d}f_{l}(x_{l};\B{\alpha}_{l})\nonumber\\
&&\qquad{}\cdot \prod_{j=1}^{d-1}\prod_{i=1}^{d-j}c_{i,i+j|v_{ij}}
\bigl(g_{i,i+j}^{(1)}(F_{i}(x_{i};\B{\alpha}_{i}),\ldots
,F_{i+j-1}(x_{i+j};\B{\alpha}_{i+j-1});\B{\theta}_{i\rightarrow
i+j-1}),\\
 &&\hphantom{\qquad{}\cdot \prod_{j=1}^{d-1}\prod_{i=1}^{d-j}c_{i,i+j|v_{ij}}
\bigl(}g_{i,i+j}^{(2)}(F_{i+1}(x_{i+1};\B{\alpha}_{i+1}),
\ldots,F_{i+j}(x_{i+j};\B{\alpha}_{i+j});\B{\theta
}_{i+1\rightarrow i+j}); \B{\theta}_{i,i+j|v_{ij}} \bigr).\nonumber
\end{eqnarray}

Recall that IFM estimates $\hat{\B{\theta}}{}^{\mathrm{IFM}}$ are obtained by plugging
the estimated marginal parameters $\hat{\B{\alpha}}{}^{\mathrm{IFM}}$ into the function
$l_{C}$. Semiparametric estimation consists in replacing the parametric marginal
c.d.f.s $u_{j}=F_{j}(x_{j};\B{\alpha}_{j})$ in $l_ {C}$ with the corresponding
empirical ones\looseness=-1
\begin{eqnarray*}
u_{jn}&=&F_{jn}(x_{j})=\frac{1}{n+1}\sum_{k=1}^{n}I(x_{jk}\leq x_{j}),\\
I(A) &=&
\cases{
1,  &\quad if \textit{A} is true,\cr
0, &\quad otherwise.
}
\end{eqnarray*}\looseness=0
The resulting pseudo log-likelihood function
$l_{C,P}(\B{\theta};\mathbf{x})$, given by
\begin{eqnarray*}
 &&l_{C,P}(\B{\theta};\mathbf{x})\\
  &&\quad= \sum_{k=1}^{n}\log(c_{1\ldots d}
(F_{1n}(x_{1k}),\ldots,F_{dn}(x_{dk});\B{\theta} ) )\\
   &&\quad=\sum_{k=1}^{n}\sum_{j=1}^{d-1}\sum_{i=1}^{d-j}\log
\bigl(c_{i,i+j|v_{ij}} \bigl(g_{i,i+j}^{(1)}(F_{i,n}(x_{ik}),\ldots
,F_{i+j-1,n}(x_{i+j-1,k});
\B{\theta}_{i\rightarrow i+j-1}), \\
 &&\hphantom{\quad=\sum_{k=1}^{n}\sum_{j=1}^{d-1}\sum_{i=1}^{d-j}\log
\bigl(c_{i,i+j|v_{ij}} \bigl(}g_{i,i+j}^{(2)}(F_{i+1,n}(x_{i+1,k}),\ldots,F_{i+j,n}(x_{i+j,k});
\B{\theta}_{i+1\rightarrow i+j}); \B{\theta}_{i,i+j|v_{ij}} \bigr) \bigr),
\end{eqnarray*}
is just a function of $\B{\theta}$. To obtain the semiparametric
estimator $\hat{\B{\theta}}{}^{\mathrm{SP}}$, one simply maximises
$l_{C,P}(\B{\theta};\mathbf{X})$ with respect to $\B{\theta}$.

Returning to the four-dimensional Student's \textit{t}-vine of Section
\ref{subsubsec:ifm}, SP estimation of this model requires
a preliminary computation of the so-called pseudo-observations
$u_{ik,n}=F_{in}(x_{ik})$, $i=1,2,3,4$, $k=1,\ldots,n$. The estimate
$\hat{\B{\theta}}{}^{\mathrm{SP}}$ is obtained by maximising
$l_{C,P}(\mathbf{x};\B{\theta})$, in this case
\begin{eqnarray*}
&&\sum_{k=1}^{n} \bigl(\log(c_{12} (u_{1k,n},u_{2k,n};\B{\theta}_{12} )
)+\cdots\\
&&\hphantom{\sum_{k=1}^{n} (}{}  +\log\bigl(c_{14|23} \bigl(g_{14}^{(1)}(u_{1k,n},\ldots,u_{3k,n};\B
{\theta}_{1\rightarrow3}),g_{14}^{(2)}(u_{2k,n},\ldots,u_{4k,n};\B
{\theta}_{2\rightarrow4});\B{\theta}_{14|23} \bigr)\bigr )\bigr ),
\end{eqnarray*}
with
$\B{\theta}_{1\rightarrow3}=(\B{\theta}_{12},\B{\theta}_{23},\B
{\theta}_{13|2})$
and
$\B{\theta}_{2\rightarrow4}=(\B{\theta}_{23},\B{\theta}_{34},\B
{\theta}_{24|3})$,
over $\B{\theta}$.

In addition to the assumptions made for the ML estimator, assume that
$c_{1\ldots d}$ fulfills condition~(A.1) from Tsukahara~\cite{tsukahara05}. Then,
the procedure corresponds to solving
\[
\frac{1}{n}\sum_{k=1}^{n}\B{\phi}^{\mathrm{SP}} (F_{1n}(X_{1k}),\ldots
, F_{1n}(X_{dk});\hat{\B{\theta}}{}^{\mathrm{SP}} )=\B{0},
\]
 with
%
\begin{equation}\label{eqn:spesteqn}
\B{\phi}_{(j-1) (d-j/2 )+i}^{\mathrm{SP}}(u_{1},\ldots,u_{d};\B
{\theta})  = \frac{\partial\log(c_{1\ldots d}(u_{1},\ldots
,u_{d};\B{\theta}) )}{\partial\B{\theta}_{i,i+j|v_{ij}}}
\end{equation}
for $i = 1,\ldots,d-j$, $j = 1,\ldots,d-1$.\vadjust{\goodbreak}

Let $\mathbf{U}$ be a \textit{d}-variate
stochastic vector distributed according to the copula $C_{1\ldots
d}(u_{1},\ldots,\break u_{d};\B{\theta})$, and define
\[
\mathbf{W}_{j}^{\mathrm{SP}}(\mathbf{U};\B{\theta}) = \int\frac{\partial^{2} \log
c_{1\ldots d}(u_{1},\ldots,u_{d};\B{\theta})}{\partial\B{\theta
}\,\partial u_{j}}I(U_{j} \leq u_{j})\,\mathrm{d}C_{1\ldots d}(u_{1},\ldots
,u_{d};\B{\theta}).
\]
Further, define
\begin{eqnarray*}
\mathbf{B}_{\theta}^{\mathrm{SP}} & =& \Var\Biggl(\sum_{j=1}^{d}\mathbf{W}_{j}^{\mathrm{SP}}(\mathbf
{U};\B{\theta}) \Biggr)+\sum_{j=1}^{d}\Cov(\B{\phi}^{\mathrm{SP}}(\mathbf{U};\B
{\theta}),\mathbf{W}_{j}^{\mathrm{SP}}(\mathbf{U};\B{\theta}) ) \\
 &=& \Var\Biggl(\sum
_{j=1}^{d}\mathbf{W}_{j}^{\mathrm{SP}}(\mathbf{U};\B{\theta})\Biggr )+\sum_{j=1}^{d}\Cov
\biggl(\frac{\partial\log c_{1\ldots d}(\mathbf{U};\B{\theta})}{\partial\B
{\theta}},\mathbf{W}_{j}^{\mathrm{SP}}(\mathbf{U};\B{\theta}) \biggr),
\end{eqnarray*}
where $\mathbf{A}=\Cov(\mathbf{Y},\mathbf{Z})$, for two stochastic vectors $\mathbf{Y}$
and $\mathbf{Z}$, is the matrix with elements
$A_{ij} = \Cov(Y_{i},Z_{j})+\Cov(Y_{j},Z_{i})$.
The matrix $\mathbf{B}_{\theta}^{\mathrm{SP}}$ quantifies the effect of replacing
the parametric marginal c.d.f.s with empirical ones. According to
Genest \textit{et~al.}~\cite{genest95} and later shown by Tsukahara~\cite{tsukahara05},
$\hat{\B{\theta}}{}^{\mathrm{SP}}$ is, under the mentioned conditions, consistent
and asymptotically normal:
%
\begin{eqnarray}\label{eqn:covsp}
 \sqrt{n} (\hat{\B{\theta}}{}^{\mathrm{SP}}-\B{\theta} ) &\stackrel
{d}{\longrightarrow} &\N(\B{0},\mathbf{V}^{\mathrm{SP}} ),\nonumber
\\[-8pt]
\\[-8pt]
 \mathbf{V}^{\mathrm{SP}}
&=& \B{\I}_{\theta}^{-1}+\B{\I}_{\theta}^{-1}\mathbf{B}_{\theta
}^{\mathrm{SP}}\B{\I}_{\theta}^{-1}.\nonumber
\end{eqnarray}
Due to the completely separate and independent estimation of
marginal and dependence parameters, the semiparametric estimator
is more robust to misspecification of the margins than ML and IFM
(Kim \textit{et~al.}~\cite{silvapulle07}). If either of the latter two produce estimates
that are
rather different from the former, it indicates that the chosen margins
or copulae are inadequate for the data.

Computationally, SP is comparable to IFM. Hence, for high-dimensional
$\B{\theta}$, although faster than ML, this procedure will require
good start values, and may still be too demanding for PCCs.
\section{PCC parameter estimators}\label{sec:pccest}
If the number of PCC parameters $\B{\theta}$ is high enough,
the estimators considered so far will be computationally too
heavy. In any case, they necessitate appropriate start values.
The next estimator, designed for pair-copula constructions,
addresses this particular issue.
\subsection{Stepwise semiparametric estimator (SSP)}\label{subsec:ssp}
As in semiparametric estimation, the marginal parameters
are handled separately, and the parametric
margins in the PCC log-likelihood function $l_{C}$
are replaced with the nonparametric ones. The idea
is to estimate the PCC parameters level by level,
conditioning on the parameters from preceding levels of
the structure. Define
\[
l_{C,P,j}(\B{\theta}_{1},\ldots,\B{\theta}_{j};\mathbf{x})  = \sum
_{k=1}^{n}\sum_{l=1}^{j}\psi_{l}(F_{1n}(x_{1k}),\ldots
,F_{dn}(x_{dk});\B{\theta}_{1},\ldots,\B{\theta}_{l}),
\]
with
%
\begin{eqnarray}\label{eqn:psi-funcD}
&&\psi_{j}(u_{1},\ldots,u_{d};\B{\theta}_{1},\ldots,\B{\theta
}_{j})\nonumber\\
&&\quad= \sum_{i=1}^{d-j}\log
\bigl(c_{i,i+j|v_{ij}}\bigl(g_{i,i+j}^{(1)}(u_{i},\ldots,u_{i+j-1};\B{\theta
}_{i\rightarrow i+j-1}),\\
& &\hphantom{\quad=\sum_{i=1}^{d-j}\log
\bigl(c_{i,i+j|v_{ij}}\bigl(} g_{i,i+j}^{(2)}(u_{i+1},\ldots,u_{i+j};\B
{\theta}_{i+1\rightarrow i+j});\B{\theta}_{i,i+j|v_{ij}}\bigr)\bigr )\nonumber
\end{eqnarray}
for $j=1,\ldots,d-1$. Hence, $l_{C,P,j}$ is the sum over all log
pair-copula densities up to, and including, level $j$. To obtain the
parameter estimates $\hat{\B{\theta}}{}^{\mathrm{SSP}}_{j}$ for level $j$, one
plugs the estimates
$\hat{\B{\theta}}{}^{\mathrm{SSP}}_{1},\ldots,\hat{\B{\theta}}{}^{\mathrm{SSP}}_{j-1}$
from preceding levels into $l_{C,P,j}$ and maximises it with respect
to $\B{\theta}_{j}$. Assuming the standard conditions for the ML
estimator are fulfilled (see Section~\ref{sec:multicopest}), this
corresponds to solving the estimating equations
$\frac{1}{n}\sum_{k=1}^{n}\B{\phi}^{\mathrm{SSP}} (F_{1n}(X_{1k}),\ldots
,F_{dn}(X_{dk});\hat{\B{\theta}}{}^{\mathrm{SSP}} )=\B{0}$,
with
%
\begin{eqnarray}\label{eqn:sspesteqn}
\B{\phi}_{(j-1) (d-j/2 )+i}^{\mathrm{SSP}}(u_{1},\ldots,u_{d};\B
{\theta}_{1},\ldots,\B{\theta}_{j}) & =& \frac{\partial}{\partial
\B{\theta}_{i,i+j|v_{ij}}}\sum_{l=1}^{j}\psi_{l}(u_{1},\ldots
,u_{d};\B{\theta}_{1},\ldots,\B{\theta}_{l})\nonumber
\\[-8pt]
\\[-8pt]
 & =& \frac
{\partial}{\partial\B{\theta}_{i,i+j|v_{ij}}}\psi_{j}(u_{1},\ldots
,u_{d};\B{\theta}_{1},\ldots,\B{\theta}_{j})\nonumber
\end{eqnarray}
for $i = 1,\ldots,d-j$, $j = 1,\ldots,d-1$. Compared to the SP
equations \eqref{eqn:spesteqn}, the full log copula density $\log
c_{1\ldots d}$ is now replaced by the sum of log copula densities up
to, and including, the level the parameter belongs to. The
corresponding estimation procedure is presented in Algorithm 1
(\ref{sup:SSPalg} of Hob{\ae}k Haff~\cite{ho2012}). If none of the
pair-copulae constituting the structure share parameters, which will
usually be the case, the estimating equations are reduced to
$\frac{\partial}{\partial\B{\theta}_{i,i+j|v_{ij}}}\log
(c_{i,i+j|v_{ij}})$. This means that the optimisation is performed for
each copula, individually.

Let us return to the four-dimensional D-vine considered in
Section~\ref{subsubsec:sp}. As in the SP procedure, one
computes the pseudo-observations $u_{ik,n}=F_{in}(x_{ik})$,
$i=1,2,3,4$, $k=1,\ldots,n$. One starts with the level 1
parameters, estimating each of the pairs $(\rho_{i,i+1},\nu_{i,i+1})$
by optimising
$\sum_{k=1}^{n}\log(c_{i,i+1}(u_{ik,n},u_{i+1,k,n};\rho_{i,i+1},\nu
_{i,i+1}))$,
for $i=1,2,3$. One subsequently computes the copula arguments
for level 2,
$u_{i|i+1,k,n}=h_{i,i+1}(u_{ik,n},u_{i+1,k,n}; \hat{\rho
}_{i,i+1},\hat{\nu}_{i,i+1}))$
and $u_{i+2|i+1,k,n}=h_{i+2,i+1}(u_{i+2,k,n},u_{i+1,k,n};\hat{\rho
}_{i+1,i+2},\break\hat{\nu}_{i+1,i+2}))$,
$i=1,2,3$, $k=1,\ldots,n$, by plugging the resulting
estimates into the adequate $h$-functions \eqref{eqn:h-func}.
At level 2, one estimates each of the pairs
$(\rho_{i,i+2|i+1},\nu_{i,i+2|i+1})$, for $i=1,2$, by
maximising
$\sum_{k=1}^{n}\log(c_{i,i+2|i+1}(u_{i|i+1,k,n}, u_{i+2|i+1,k,n};\rho
_{i,i+2|i+1},\nu_{i,i+2|i+1}))$.
Next, one computes the copula arguments $u_{1|23,k,n}$
and $u_{4|23,k,n}$ for level 3 by plugging the estimates from
level 2 into $h_{13|2}$ and $h_{24|3}$. Finally, to obtain
$(\hat{\rho}_{14|23},\hat{\nu}_{14|23})$, one optimises $\sum
_{k=1}^{n}\log(c_{14|23}(u_{1|23,k,n},u_{4|23,k,n};\rho_{14|23},\nu
_{14|23}))$.

When some of the copulae share parameters, the procedure is
a little different. Let us for instance consider a four-dimensional
Student's \textit{t}-copula with correlations
$(\rho_{12},\rho_{23},\rho_{34}, \rho_{13},\rho_{24},\rho_{14})$
and $\nu$ degrees of freedom. This is also a D-vine consisting of
Student's \textit{t}-copulae (see for instance Min and Czado~\cite{minczado10}). The
correlation parameters of these copulae are now the corresponding
partial correlations
$(\rho_{12},\rho_{23},\rho_{34},\rho_{13|2},\rho_{24|3},\rho_{14|23})$.
However, the degrees of freedom parameter is shared. More
specifically, it is $\nu$ for the three copulae at the
ground level, $\nu+1$ at level 2 and $\nu+2$
for the top level copula. The SSP estimation procedure
is now as follows. Having computed the pseudo-observations,
one maximises the level 1 function
\[
\sum_{k=1}^{n}\psi_{1}(u_{1k,n},\ldots,u_{4k,n};\rho_{12},\rho
_{23},\rho_{34},\nu) = \sum_{k=1}^{n}\sum_{i=1}^{3} \log
(c_{i,i+1}(u_{i,k,n},u_{i+1,k,n};\rho_{i,i+1},\nu) )
\]
over $(\rho_{12},\rho_{23},\rho_{34},\nu)$. Then, one calculates
the copula arguments for level 2 as described above. At the second
level, one estimates $\rho_{13}$ and $\rho_{24}$, which are not shared
by $c_{13|2}$ and $c_{24|3}$. More specifically, one optimises each of
$\sum_{k=1}^{n}\log(c_{i,i+2|i+1}(u_{i|i+1,k,n},u_{i+2|i+1,k,n};\break
\rho_{i,i+2},\hat{\rho}_{i,i+1},\hat{\rho}_{i+1,i+2},\hat{\nu}))$,
over $\rho_{i,i+2}$, $i=1,2$ (note that
$\hat{\rho}_{i,i+1},\hat{\rho}_{i+1,i+2}$ are needed to compute
the partial correlations $\rho_{i,i+2|i+1}$). Next, one computes
the copula arguments for the top level copula, and finally, one
maximises
$\sum_{k=1}^{n}\log(c_{14|23}(u_{1|23,k,n},u_{4|23\},k,n};\rho
_{14},\hat{\rho}_{12},\break \ldots,\hat{\rho}_{24|3},\hat{\nu}))$
over $\rho_{14}$.
Note however that although it is possible to estimate the parameters
of a multivariate Student's \textit{t}-copula as described above, it is
unnecessarily complex. In practice, one would typically estimate
the correlation parameters via the corresponding Kendall's $\tau$
coefficients, and subsequently optimise the pseudo
log-likelihood function $l_{C,P}$ over $\nu$, plugging in the
estimated correlations, as described in for instance
McNeil \textit{et~al.}~\cite{mcneil06} (page 231). The main purpose of the PCC is to model pairs
with different behaviour. If one does not really need that
flexibility, then using a PCC is like using a sledgehammer
to crack a nut.

Let us now consider conditions (A.1)--(A.5) from Tsukahara~\cite{tsukahara05}.
The last four of these are the standard conditions for the ML
estimator, but on the score functions \eqref{eqn:sspesteqn}. Further,
define
\begin{eqnarray*}
\B{\phi}_{(j-1) (d-j/2 )+i}(\mathbf{u};\B{\theta}_{1},\ldots
,\B{\theta}_{j}) &=& \frac{\partial}{\partial\B{\theta
}_{i,i+j|v_{ij}}}\psi_{j}(\mathbf{u};\B{\theta}_{1},\ldots,\B{\theta
}_{j}) \equiv\B{\psi}_{ij,\theta}(\mathbf{u};\B{\theta}_{1},\ldots
,\B{\theta}_{j}), \\
\frac{\partial}{\partial u_{k}}\B{\phi}_{(j-1)
(d-j/2 )+i}(\mathbf{u};\B{\theta}_{1},\ldots,\B{\theta}_{j})
&=& \frac{\partial}{\partial u_{k}}\B{\psi}_{ij,\theta}(\mathbf{u};\B
{\theta}_{1},\ldots,\B{\theta}_{j}) \equiv\B{\psi}_{ij,\theta
,u_{k}}(\mathbf{u};\B{\theta}_{1},\ldots,\B{\theta}_{j}).
\end{eqnarray*}
Let $\Q$ and $\Rr$ be the sets of positive, symmetric, inverse square
integrable functions on $[0,1]$ and reproducing u-shaped functions
on $[0,1]$, respectively, as defined in Tsukahara~\cite{tsukahara05}. Further,
let $|\B{\theta}_{ij|v_{ij}}|=l_{ij}$ be the number of parameters of
the pair-copula $C_{i,i+j|v_{ij}}$. For the SSP estimator, Condition
(A.1) may then be phrased in the following way (note that a
subscript `\textit{j}' on $\phi$ is missing in Tsukahara~\cite{tsukahara05}).
\begin{Cond}\label{cond:ssp}
For each $\B{\theta}$,
$\B{\psi}_{ij,\theta}=(\psi_{ij,\theta,1},\ldots,
\psi_{ij,\theta
,l_{ij}})$
and $\B{\psi}_{ij,\theta,u_{k}}=(\psi_{ij,\theta,u_{k},1},
\ldots,\break\psi_{ij,\theta,u_{k},l_{ij}}), j=1,\ldots,d-1, i=1,\ldots,d-j$ are
continuous, and there exist functions
$r_{ij,k},\tilde{r}_{ij,k} \in\Rr$ and $q_{ij,k} \in\Q$,
such that
\begin{eqnarray*}
|\psi_{ij,\theta,m}(\mathbf{u};\B{\theta}_{1},\ldots,\B{\theta
}_{j})| & \leq&\prod_{l=1}^{d}r_{ij,l}(u_{l}),\\
 |\psi_{ij,\theta
,u_{k},m}(\mathbf{u};\B{\theta}_{1},\ldots,\B{\theta}_{j})| & \leq&
\tilde{r}_{ij,k}(u_{k})\prod_{l\neq k}r_{ij,l}(u_{l})
\end{eqnarray*}
for $k,l=1,\ldots,d$, $j=1,\ldots,d-1$, $i=1,\ldots,d-j$,
$m=1,\ldots,l_{ij}$, with
\begin{eqnarray*}
\int\Biggl(\prod_{l=1}^{d}r_{ij,l}(u_{l})\Biggr )^{2}\,\mathrm{d}C_{1\ldots d}(u_{1},\ldots
,u_{d};\B{\theta}) & <& \infty, \\
\int\biggl(q_{ij,k}(u_{k})\tilde
{r}_{ij,k}(u_{k})\prod_{l\neq k}r_{ij,l}(u_{l}) \biggr)^{2}\,\mathrm{d}C_{1\ldots
d}(u_{1},\ldots,u_{d};\B{\theta}) & <& \infty.
\end{eqnarray*}
\end{Cond}

When none of the pair-copulae share parameters, Condition
\ref{cond:ssp} becomes a condition on each of them, individually.

Once more let $\mathbf{U}$ be distributed according to
$C_{1\ldots d}(u_{1},\ldots,u_{d};\B{\theta})$, as well as
$\B{\psi}_{\theta}=(\B{\psi}_{11,\theta},\ldots,\B{\psi
}_{1,d-1,\theta})$
and
$\B{\psi}_{\theta,u_{j}}=(\B{\psi}_{11,\theta,u_{j}},\ldots,\B
{\psi}_{1,d-1,\theta,u_{j}})$.
Define
\begin{eqnarray*}
\mathbf{W}_{j}^{\mathrm{SSP}}(\mathbf{U};\B{\theta}) & =& \int\frac{\partial
}{\partial u_{j}}\B{\phi}^{\mathrm{SSP}}(u_{1},\ldots,u_{d};\B{\theta
})I(U_{j} \leq u_{j})\,\mathrm{d}C_{1\ldots d}(u_{1},\ldots,u_{d};\B{\theta}) \\
&= &\int\B{\psi}_{\theta,u_{j}}(u_{1},\ldots,u_{d};\B{\theta
})I(U_{j} \leq u_{j})\,\mathrm{d}C_{1\ldots d}(u_{1},\ldots,u_{d};\B{\theta})
\end{eqnarray*}
and the matrix
\begin{eqnarray*}
\mathbf{B}_{\theta}^{\mathrm{SSP}} & =& \Var\Biggl(\sum_{j=1}^{d}\mathbf{W}_{j}^{\mathrm{SSP}}(\mathbf
{U};\B{\theta}) \Biggr)+\sum_{j=1}^{d}\Cov(\B{\phi}^{\mathrm{SSP}}(\mathbf{U};\B
{\theta}),\mathbf{W}_{j}^{\mathrm{SSP}}(\mathbf{U};\B{\theta}) )\\
& =& \Var\Biggl(\sum_{j=1}^{d}\mathbf{W}_{j}^{\mathrm{SSP}}(\mathbf{U};\B{\theta}) \Biggr)+\sum
_{j=1}^{d}\Cov(\B{\psi}_{\theta}(\mathbf{U};\B{\theta}),\mathbf
{W}_{j}^{\mathrm{SSP}}(\mathbf{U};\B{\theta}) ).
\end{eqnarray*}
Moreover, define the two matrices
\begin{eqnarray*}
\B{\K}_{\theta} & =& \E(\B{\phi}^{\mathrm{SSP}} (\B{\phi}^{\mathrm{SSP}} )^{T} ) =
\pmatrix{
\B{\K}_{\theta,1,1} & & & \B{0}\cr
 & \ddots& & \cr
 \B{0}^{T} & & \B{\K}_{\theta,d-2,d-2} & \B{0}\cr
  \B{0}^{T} & & \B{0}^{T} & \B{\I}_{\theta,d-1,d-1}},\\
\B{\J}_{\theta} & =& \E\biggl(-\frac{\partial\B{\phi
}^{\mathrm{SSP}}}{\partial\B{\theta}^{T}} \biggr) =
\pmatrix{
\B{\J}_{\theta,1,1} & & & \B{0} \cr
\vdots& & & \cr
\vdots& \ddots& &\cr
 \B{\J}_{\theta,d-2,1} & \cdots& \B{\J}_{\theta,d-2,d-2} & \B{0}\cr
 \B{\I}_{\theta,d-1,1} & \cdots& \B{\I}_{\theta,d-1,d-2} & \B{\I}_{\theta,d-1,d-1}},
\end{eqnarray*}
where the blocks
$\B{\K}_{\theta,i,j}=\E( (\frac{\partial\psi_{i}}{\partial\B
{\theta}_{i}} ) (\frac{\partial\psi_{j}}{\partial\B{\theta}_
{j}} )^{T} )$
and
$\B{\J}_{\theta,i,j}=-\E(\frac{\partial^{2} \psi_{i}}{\partial
\B{\theta}_{i}\,\partial\B{\theta}_ {j}^{T}} )$,
$i,j=1,\ldots,d-1$, correspond to each of the construction's levels.
The block diagonal and block lower triangular forms of $\B{\K
}_{\theta}$
and $\B{\J}_{\theta}$, respectively, follow from the structure of
the estimating equations (see Appendix~\ref{subsec:sspcomp}).
More specifically, the $\psi$ functions depend on all the parameters
from previous levels but not from following levels. Further,
the estimating equations for the top level copula parameters are
based on the full copula, as for the SP estimator. This accounts
for the appearance of blocks from the Fisher matrix $\B{\I}_{\theta}$
in the last rows of $\B{\K}_{\theta}$ and $\B{\J}_{\theta}$. If all
pair-copulae are from one-parameter families, then $\B{\K}_{\theta}$
and $\B{\J}_{\theta}$ are $d(d-1)/2 \times d(d-1)/2$ matrices.

We now have all the necessary components to establish the
asymptotic properties of the stepwise semiparametric estimator.
\begin{Theorem}\label{theo:ssp}
Under Condition~\ref{cond:ssp}, as well as Conditions
\textup{(A.2)--(A.5)} of Tsukahara~\cite{tsukahara05}, the SSP estimator
$\hat{\B{\theta}}{}^{\mathrm{SSP}}$ is consistent for $\B{\theta}$ and
asymptotically normal:
%
\begin{eqnarray}\label{eqn:covssp}
 \sqrt{n} (\hat{\B{\theta}}{}^{\mathrm{SSP}}-\B{\theta} )& \stackrel
{d}{\longrightarrow}& \N(\B{0},\mathbf{V}^{\mathrm{SSP}} ),\nonumber
\\[-8pt]
\\[-8pt]
 \mathbf{V}^{\mathrm{SSP}}& =& \B{\J}_{\theta}^{-1}\B{\K}_{\theta} (\B{\J}_{\theta
}^{-1} )^{T}+\B{\J}_{\theta}^{-1}\mathbf{B}_{\theta}^{\mathrm{SSP}} (\B{\J
}_{\theta}^{-1} )^{T}.\nonumber
\end{eqnarray}
\end{Theorem}
\begin{pf}
Theorem~\ref{theo:ssp} follows directly Theorem 1 of Tsukahara~\cite{tsukahara05},
with the estimating equations~\eqref{eqn:sspesteqn}. Note that Theorem 1
of Tsukahara~\cite{tsukahara05} is valid for the multiparameter case $m > 1$, despite
some misprints and imprecisions in the original paper. Specifically,
Condition (A.1) is assumed valid for every element of the vector $\B
{\phi}$
of estimating equations, that is, a subscript `\textit{j}' is needed.
Proposition 3
holds for $m > 1$ thanks to the Cramer--Wold device. The rest of the argument
works for $m > 1$, using $\Vert\cdot\Vert$ instead of $|\cdot|$ for the norm
(Tsukahara~\cite{tsukahara2011}).
\end{pf}

In order to construct confidence intervals for $\B{\theta}$,
one needs a consistent estimate of $\mathbf{V}^{\mathrm{SSP}}$. As noted in
Tsukahara~\cite{tsukahara05}, one may estimate this covariance matrix
consistently by replacing expectations and variances in
\eqref{eqn:covssp} by sample equivalents, and plugging in the
estimate $\hat{\B{\theta}}{}^{\mathrm{SSP}}$. More specifically, letting
$\B{\psi}_{\theta,\theta} = \frac{\partial}{\partial\B{\theta
}^{T}}\B{\psi}_{\theta}\mathbf{W}^{\mathrm{SSP}}=\sum_{j=1}^{d}\mathbf
{W}_{j}^{\mathrm{SSP}}$, and
$u_{ik}$, $i=1,\ldots,d$, $k=1,\ldots,n$, are the
pseudo-observations,
\[
\hat{\mathbf{V}}^{\mathrm{SSP}}  = \hat{\B{\J}}{}^{-1}_{\theta}
\hat{\B{\K}}_{\theta} (\hat{\B{\J}}{}^{-1}_{\theta} )^{T}+
\hat{\B{\J}}{}^{-1}_{\theta}\hat{\mathbf{B}}{}^{\mathrm{SSP}}_{\theta}
(\hat{\B{\J}}{}^{-1}_{\theta} )^{T},
\]
with
\begin{eqnarray*}
\hat{\B{\K}}_{\theta} & =& \frac{1}{n}\sum_{k=1}^{n}\B{\psi
}_{\theta} (u_{1k},\ldots,u_{dk};\hat{\B{\theta}}{}^{\mathrm{SSP}} ) (\B
{\psi}_{\theta} (u_{1k},\ldots,u_{dk};\hat{\B{\theta}}{}^{\mathrm{SSP}} )
)^{T},\\
 \hat{\B{\J}}_{\theta} & =& -\frac{1}{n}\sum_{k=1}^{n}\B
{\psi}_{\theta,\theta} (u_{1k},\ldots,u_{dk};\hat{\B{\theta
}}^{\mathrm{SSP}} ),\\[-2pt]
\hat{\mathbf{B}}{}^{\mathrm{SSP}}_{\theta} & =& \frac{1}{n}\sum
_{k=1}^{n}\mathbf{W}^{\mathrm{SSP}} (u_{1k},\ldots,u_{dk};\hat{\B{\theta}}{}^{\mathrm{SSP}}
) (\mathbf{W}^{\mathrm{SSP}} (u_{1k},\ldots,u_{dk};\hat{\B{\theta}}{}^{\mathrm{SSP}} )
)^{T} \\[-2pt]
&&{} +\frac{1}{n}\sum_{k=1}^{n}\sum_{j=1}^{d}\B{\psi}_{\theta}
(u_{1k},\ldots,u_{dk};\hat{\B{\theta}}{}^{\mathrm{SSP}} ) (\mathbf{W}^{\mathrm{SSP}}
(u_{1k},\ldots,u_{dk};\hat{\B{\theta}}{}^{\mathrm{SSP}} ) )^{T}.
\end{eqnarray*}
In most cases, there is no analytic
expression for the derivatives $\B{\psi}_{\theta}$ and
$\B{\psi}_{\theta,\theta}$, but they can be approximated
numerically. However, the computation of the
$\mathbf{W}_{j}$-vectors involves $d$-dimensional integrals, which
is more problematic. In practice, one will not be able to
compute the above covariance matrix estimate for $d > 3$.
Instead, one will have to resort to some resampling technique,
such as parametric bootstrap from
$C_{1\ldots d}(\cdot;\hat{\bolds{\theta}}{}^{\mathrm{SSP}})$, as described in Example
\ref{ex:raindata} of Section~\ref{sec:examples}.\vspace*{-2pt}

\begin{Theorem}\label{theo:sspnorm}
Under the conditions of Theorem~\ref{theo:ssp}, the SSP estimator
$\hat{\B{\theta}}{}^{\mathrm{SSP}}$ is asymptotically semiparametrically
efficient for the parameters $\B{\theta}$ of the Gaussian
copula.\vspace*{-2pt}
\end{Theorem}

The proof is given in Appendix~\ref{subsec:proofsspnorm}.

In general, the stepwise semiparametric estimator
$\hat{\B{\theta}}{}^{\mathrm{SSP}}$ is asymptotically less efficient than
$\hat{\bolds{\theta}}{}^{\mathrm{SP}}$, since it at a given level discards all
information from following levels. Nonetheless, the levelwise
estimation significantly improves the computational efficiency.
The SSP estimator is therefore adequate for medium to
high-dimensional models, and to produce start values for the SP
estimator. Further, a substantial difference between SSP and
SP estimates may be a sign that the copulae are unsuitable.
Hence, one may use the SSP estimator to assess the sensitivity
to the chosen copulae. Moreover, it is inherently suited for
determining an appropriate PCC for a data set, which consists
in choosing an ordering of the variables and a set of parametric
pair-copulae in a stepwise manner. Once the ordering is fixed,
one finds suitable copulae for the ground level, based on the
pseudo-observations. At the second level, the necessary pair-copula
arguments are obtained by transforming the pseudo-observations
with the adequate $h$-functions, which depend on the chosen
ground level copulae. This requires ground level parameter
estimates, which can be provided by the SSP estimator. After
one has selected copulae for the second level, one proceeds in
the same manner for the remaining levels. Of course, one could
construct a similar, stepwise estimator with a different
transformation to uniform margins, for instance using the
parametric margins as in IFM estimation. That particular
estimator was in fact proposed by Joe and Xu
\cite{joexu96}.\vspace*{-2pt}

\subsection{Robustness}\label{subsec:robustness}
The SSP estimator is a substantial improvement over the three
former in terms of computational speed. However, it presupposes
that the specified model is the true one. If the amount of data
available is high enough, it should, in most cases, be possible
to find adequate marginal distributions. For the pair-copulae,
the task is more complex.\vadjust{\goodbreak} Using the pseudo-observations, one
may obtain a reasonable model for the ground level. Subsequently,
however, one must condition on choices from previous levels, as
described above. One would therefore expect the quality of the
model to decrease with the construction level.

SSP estimation consists in replacing the parametric margins in
the function $l_{C}$ with the nonparametric ones,
while keeping the parametric forms of the conditional distributions,
that is, the $g$-functions \eqref{eqn:g-func}. The resulting estimator
is robust toward misspecification of the margins, but not of the
pair-copulae. By replacing also the conditional distributions
with nonparametric versions, one would reduce this sensitivity
to chosen pair-copulae preceding in the structure. One possibility
is the empirical conditional distribution proposed by
Stute~\cite{stute86}:
%
\begin{equation}\label{eqn:empcondmod}
F_{i|v,n}(x_{i}|\mathbf{x}_{v}) = \frac{n}{n+1}\frac{\sum
_{k=1}^{n}I(x_{ik} \leq x_{i})K_{s} (\mathbf{x}_{v}-\mathbf{x}_{v,k} )}{\sum
_{k=1}^{n}K_{s} (\mathbf{x}_{v}-\mathbf{x}_{v,k} )}, \qquad K_{s}(\mathbf{y}) = \frac
{1}{s^{l}}K \biggl(\frac{\mathbf{y}}{s} \biggr),
\end{equation}
where $l$ is the dimension of $\mathbf{x}_{v}$, $K$ is a kernel function on
$\R^{l}$ and $s$ the bandwidth parameter. The definition
\eqref{eqn:empcondmod} is slightly modified here to avoid boundary problems
in $0$ and $1$. Provided $h \rightarrow0$ and $ns^{l} \rightarrow
\infty$ as
$n\rightarrow\infty$, it converges almost surely to the true conditional
distribution, though at a rather slow pace of order $(ns^{l})^{1/2}$. The
quality of the estimates will therefore significantly decrease with the
level number. Alternative versions of the empirical conditional distribution
function, such as the one proposed by Hall and Yao~\cite{hall05}, share this unfortunate
property.

Recall that the conditional distributions of interest are recursions
of $h$-functions \eqref{eqn:h-func}, which, in turn, are conditional
distributions of uniform variables with a conditioning set of length
one. These functions can therefore be estimated nonparametrically by
\eqref{eqn:empcondmod} with $l=1$. Seemingly, one can exploit
this to avoid the curse of dimensionality. However, the two arguments
of the $h$-functions are again $h$-functions from the preceding level.
Hence, the error propagates from level to level, and as expected, the
resulting rate of convergence is of the same order as for the
original variables, that is $(ns^{l})^{1/2}$.

Accordingly, the estimator suggested above becomes unreliable
already at the fourth or fifth level of the structure,
depending on the amount of data. Since the intention is to
improve the quality of estimates at higher levels, it is
in practice useless, unless the rate of convergence is
increased by additional assumptions on the conditional
distributions.
\subsection{C-vines}\label{subsec:c-vines}
For simplicity, we have only considered D-vines so far.
The same results are however easily obtained for C-vines
(see Figure~\ref{fig:CDvine}) and more general regular
vines (though computation of the log-likelihood
function is more complex (Dissmann \textit{et~al.}~\cite{dissmann11})).

The p.d.f. of a C-vine is given by (Aas \textit{et~al.}~\cite{Vines})
%
\begin{eqnarray}\label{eqn:cvinepdf}
&&f_{1\ldots d}(x_{1},\ldots,x_{d};\B{\alpha},\B{\theta})\nonumber\\
&&\quad =  \prod
_{l=1}^{d}f_{l}(x_{l};\B{\alpha}_{l})\nonumber
\\[-8pt]
\\[-8pt]
&&\qquad{}\cdot\prod_{j=1}^{d-1}\prod
_{i=1}^{d-j}c_{j,j+i|z_{j-1,0}} (F_{j|z_{j-1,0}}(x_{j}|\mathbf
{x}_{z_{j-1,0}};\B{\alpha}_{z_{j0}},\B{\theta}_{j\_0}),\nonumber\\[-2pt]
& &\hphantom{\qquad{}\cdot\prod_{j=1}^{d-1}\prod
_{i=1}^{d-j}c_{j,j+i|z_{j-1,0}} (}
F_{j+i|z_{j-1,0}}(x_{j+i}|\mathbf{x}_{z_{j-1,0}};\B{\alpha}_{z_{ji}},\B
{\theta}_{j\_i});\B{\theta}_{j,j+i|z_{j-1,0}} ),\nonumber
\end{eqnarray}
where $z_{ji}=\{1,\ldots,j-1,j+i\}$ and
$\B{\theta}_{j\_i}=\{\B{\theta}_{s,s+t|z_{s-1,0}}\dvt (s,s+t) \in
z_{ji}\}$,
for $0\leq i \leq d-j$, $0\leq j \leq d-1$, with
$z_{0i} = \varnothing$, $z_{1i} = \{i+1\}$ and
$\B{\theta}_{0\_i} = \B{\theta}_{1\_i} = \varnothing$. Hence, the
log-likelihood function of $n$ independent observations
from a C-vine is
%
\begin{eqnarray}\label{eqn:loglikCvine}
l(\B{\alpha},\B{\theta};\mathbf{x}) &=& \sum_{k=1}^{n}\log
(f_{1\ldots d}(x_{1k},\ldots,x_{dk};\B{\alpha},\B{\theta}))\nonumber\\[-2pt]
 &=& \sum_{k=1}^{n}\sum_{l=1}^{d}\log(f_{l}(x_{lk};\B{\alpha
}_{l}))\nonumber\\[-2pt]
&&{}+ \sum_{k=1}^{n}\sum_{j=1}^{d-1}\sum
_{i=1}^{d-j}\log(c_{j,j+i|z_{j-1,0}} (F_{j|z_{j-1,0}}(x_{j}|\mathbf
{x}_{z_{j-1,0}};\B{\alpha}_{z_{j0}},\B{\theta}_{j\_0}), \\[-2pt]
&&\hphantom{{}+ \sum_{k=1}^{n}\sum_{j=1}^{d-1}\sum
_{i=1}^{d-j}\log(c_{j,j+i|z_{j-1,0}} (} F_{j+i|z_{j-1,0}}(x_{j+i}|\mathbf{x}_{z_{j-1,0}};\B{\alpha}_{z_{ji}},\B
{\theta}_{j\_i});\B{\theta}_{j,j+i|z_{j-1,0}} ) )\nonumber \\[-2pt]
& =&
l_{M}(\B{\alpha};\mathbf{x})+l_{C}(\B{\alpha},\B{\theta};\mathbf{x}).\nonumber
\end{eqnarray}
Replacing $l_{C}$ from \eqref{eqn:loglik} with $l_{C}$ from
\eqref{eqn:loglikCvine}, one retrieves the results
from Section~\ref{sec:multicopest} for C-vines. To achieve the SSP
estimator, one must simply replace the psi-function
\eqref{eqn:psi-funcD} in the estimating equations \eqref{eqn:sspesteqn}
with
%
\begin{eqnarray}\label{eqn:psi-funcC}
&&\psi_{j}(u_{1},\ldots,u_{d};\B{\theta}_{1},\ldots,\B{\theta
}_{j})\nonumber\\[-2pt]
& &\quad=  \sum_{i=1}^{d-j}\log(c_{j,j+i|z_{j-1,0}}
(F_{j|z_{j-1,0}}(x_{j}|\mathbf{x}_{z_{j-1,0}};\B{\alpha}_{z_{j0}},\B
{\theta}_{j\_0}), \\[-2pt]
&&\hphantom{\quad=  \sum_{i=1}^{d-j}\log(c_{j,j+i|z_{j-1,0}}
(} F_{j+i|z_{j-1,0}}(x_{j+i}|\mathbf{x}_{z_{j-1,0}};\B
{\alpha}_{z_{ji}},\B{\theta}_{j\_i});\B{\theta}_{j,j+i|z_{j-1,0}}
) ).\nonumber
\end{eqnarray}
Also, the $h$-functions \eqref{eqn:h-func} are redefined as
%
\begin{equation}\label{eqn:h-funccvine}
h_{j+i,j|z_{j-1,0}}(u_{j+i},u_{j}) \equiv\frac{\partial
C_{j,j+i|z_{j-1,0}}(u_{j},u_{j+i})}{\partial u_{j}}
\end{equation}
for $i=1,\ldots,d-j$, $j=1,\ldots,d-1$. The estimation procedure for a
C-vine is described in Algorithm 2 (\ref{sup:SSPalg} of Hob{\ae}k Haff
\cite{ho2012}).
%
\section{Examples}\label{sec:examples}
To compare the four estimators' performance, we have
carried out asymptotic computations on a few examples
(Examples~\ref{ex:normdist} to~\ref{ex:student}). We have
also fitted a D-vine to a set of precipitation series, using
each of the estimators (Example~\ref{ex:raindata}).\vadjust{\goodbreak}

\begin{Ex}\label{ex:normdist}
Consider the three-dimensional Gaussian distribution
\[
\pmatrix{
X_{1}\cr
 X_{2}\cr X_{3}}
 \sim\N_{3} (\B{0},\mathbf{S}\mathbf{R}\mathbf{S} ),\qquad\mathbf{S} =
\pmatrix{
\sigma_{1} & 0 & 0 \cr
0 & \sigma_{2} & 0 \cr
0 & 0 & \sigma_{3}}
, \qquad\mathbf{R} =
\pmatrix{
1 & \rho_{12} & \rho_{13}\cr
 \rho_{12} & 1 & \rho_{23}\cr
  \rho_{13} &
\rho_{23} & 1}.
\]
This distribution can be represented by a D-vine consisting of
Gaussian pair-copulae and margins, more specifically
\[
f_{123}(x_{1},x_{2},x_{3};\mathbf{S},\mathbf{R}) = c(u_{1},u_{2};\rho
_{12})c(u_{2},u_{3};\rho_{23})c(u_{1|2},u_{3|2};\rho_{13|2})\prod
_{i=1}^{3}f(x_{i};\sigma_{i}),
\]
where
\begin{eqnarray*}
c(u_{i},u_{j};\rho) &= & \frac{\exp\{-\rho/(2(1-\rho
^{2}))(\rho\Phi^{-1}(u_{i})^{2}+\rho\Phi^{-1}(u_{j})^{2}-2\Phi
^{-1}(u_{i})\Phi^{-1}(u_{j})) \}}{\sqrt{1-\rho^{2}}}\\
 f(x;\sigma) &= &
\frac{\exp\{-{x^{2}}/(2\sigma^{2}) \}}{\sqrt{2\uppi}\sigma},\qquad
\rho_{13|2} = \frac{\rho_{13}-\rho_{12}\rho_{23}}{\sqrt{(1-\rho
_{12}^{2})(1-\rho_{23}^{2})}},\\
 u_{i|j} &= & h_{ij}(u_{i},u_{j};\rho) =
\Phi\biggl(\frac{\Phi^{-1}(u_{i})-\rho\Phi^{-1}(u_{j})}{\sqrt{1-\rho
^{2}}} \biggr),
\end{eqnarray*}
with $u_{i} = \Phi(x_{i})$, $\Phi$ being the c.d.f. of the standard
Gaussian distribution.
Note that this is one of the three possible decompositions of
$f_{123}$.

In practice, there are scarcely any other models for which it
is feasible to do all computations analytically. It is also one
of the few distributions the IFM and SP estimators are asymptotically
efficient for, as explained below.

The ML estimators $\hat{\B{\alpha}}{}^{\mathrm{ML}}$ and $\hat{\B{\theta}}{}^{\mathrm{ML}}$
are of course the empirical standard deviations and correlations,
respectively. It is easily verified that for this particular model, the
IFM estimators $\hat{\B{\alpha}}{}^{\mathrm{IFM}}$ and $\hat{\B{\theta}}{}^{\mathrm{IFM}}$
are identical to $\hat{\B{\alpha}}{}^{\mathrm{ML}}$ and $\hat{\B{\theta}}{}^{\mathrm{ML}}$.
Thus, they are asymptotically efficient. Moreover, the SP estimator
$\hat{\B{\theta}}{}^{\mathrm{SP}}$ is semiparametrically efficient for $\B
{\theta}$,
as shown by Klaassen and Wellner~\cite{klaassen97}.

For SSP, we must compute the matrices $\B{\K}_{\theta}$, $\B{\J
}_{\theta}$
and $\mathbf{B}_{\theta}^{\mathrm{SSP}}$, defined in Section~\ref{subsec:ssp}.
The covariance matrix $\mathbf{V}^{\mathrm{SSP}}$
of $\rho_{12}, \rho_{23}$ and $\rho_{13}$, in that order
(corresponding to the PCC levels), is shown in Appendix
\ref{subsec:normexcomp}, along with $\mathbf{V}^{\mathrm{ML}}$. We
see that $\mathbf{V}^{\mathrm{SSP}}=\mathbf{V}^{\mathrm{ML}}=\mathbf{V}^{\mathrm{SP}}$. As the SP estimator is
asymptotically semiparametrically efficient for $\B{\theta}$,
so must the SSP estimator be.
\end{Ex}
\begin{Ex}\label{ex:expgumb}
Consider the three-dimensional PCC with exponential margins and
Gumbel pair-copulae:
\[
f_{123}(x_{1},x_{2},x_{3};\B{\lambda},\B{\delta}) =
c(u_{1},u_{2};\delta_{12})c(u_{2},u_{3};\delta
_{23})c(u_{1|2},u_{3|2};\delta_{13|2})\prod_{i=1}^{3}f(x_{i};\lambda_{i}),
\]
where
\begin{eqnarray*}
c(u_{i},u_{j};\delta) &= & \frac{\exp\{-(\tilde{u}_{i}^{\delta
}+\tilde{u}_{j}^{\delta})^{1/\delta}\}(\tilde{u}_{i}\tilde
{u}_{j})^{\delta-1}((\tilde{u}_{i}^{\delta}+\tilde{u}_{j}^{\delta
})^{1/\delta}+\delta-1)}{u_{i}u_{j}(\tilde{u}_{i}^{\delta}+\tilde
{u}_{j}^{\delta})^{2-1/\delta}}, \\[-2pt]
f(x;\lambda) &= & \lambda\exp\{
-\lambda x\},\qquad u_{i|j} = h_{ij}(u_{i},u_{j};\delta) = \frac{\exp\{
-(\tilde{u}_{i}^{\delta}+\tilde{u}_{j}^{\delta})^{1/\delta}(\tilde
{u}_{j})^{\delta-1}\}}{u_{j}(\tilde{u}_{i}^{\delta}+\tilde
{u}_{j}^{\delta})^{1-1/\delta}},
\end{eqnarray*}
with $u_{j} = 1-\exp\{-\lambda_{j}x_{j}\}$ and
$\tilde{u}_{j} = -\log(u_{j}), i,j = 1,2,3$. For various parameter
sets, we have computed the covariance matrices by numerical derivation
and integration. Since the dependence parameters $\B{\delta}$ are
\begin{table}
\def\arraystretch{0.9}
\tablewidth=\textwidth
\tabcolsep=0pt
  \caption{Asymptotic relative efficiencies of $\hat{\delta}_{12}$
and $\hat{\delta}_{13|2}$ from Example \protect\ref{ex:expgumb}, for various
parameter sets}\label{tab:expgumb}
\begin{tabular*}{\textwidth}{@{\extracolsep{\fill}}lllllll@{}}
  \hline
   & \multicolumn{2}{l}{\textsc{IFM}} & \multicolumn{2}{l}{\textsc{SP}}& \multicolumn{2}{l}{\textsc{SSP}} \\[-5pt]
   & \multicolumn{2}{l}{\hrulefill} & \multicolumn{2}{l}{\hrulefill}& \multicolumn{2}{l@{}}{\hrulefill}  \\
   & $\hat{\delta}_{12}$ & $\hat{\delta}_{13|2}$ & $\hat{\delta}_{12}$ & $\hat{\delta}_{13|2}$ & $\hat{\delta}_{12}$ & $\hat{\delta}_{13|2}$\\
  \hline
$(\delta_{12},\delta_{13|2})=(1.2,1.2)$ & 0.997 & 0.997 & 0.921 & 0.955 & 0.904 & 0.953\\
$(\delta_{12},\delta_{13|2})=(1.2,2)$ & 0.985 & 0.996 & 0.902 & 0.984 & 0.891 & 0.981  \\
$(\delta_{12},\delta_{13|2})=(1.2,3)$ & 0.971 & 0.994 & 0.846 & 0.990 & 0.837 & 0.987 \\
$(\delta_{12},\delta_{13|2})=(2,1.2)$ & 0.995 & 0.985 & 0.913 & 0.851 & 0.879 & 0.843\\
$(\delta_{12},\delta_{13|2})=(2,2)$ & 0.981 & 0.983 & 0.896 & 0.950 & 0.850 & 0.936 \\
$(\delta_{12},\delta_{13|2})=(2,3)$ & 0.956 & 0.969& 0.832 & 0.976 & 0.815 & 0.962\\
$(\delta_{12},\delta_{13|2})=(3,1.2)$ & 0.995 & 0.974 & 0.912 & 0.814 & 0.861 & 0.808\\
$(\delta_{12},\delta_{13|2})=(3,2)$ & 0.973 & 0.954 & 0.871 & 0.921 & 0.843 & 0.887 \\
$(\delta_{12},\delta_{13|2})=(3,3)$ & 0.944 & 0.932 & 0.825 & 0.951 & 0.777 & 0.931\\
  \hline
\end{tabular*} \vspace*{-3pt}
\end{table}
our primary interest, we let $\lambda_{1}=\lambda_{2}=\lambda_{3}=1$ in
all sets. Moreover, we let $\delta_{12}=\delta_{23}$. Table
\ref{tab:expgumb} shows the resulting asymptotic relative efficiencies
of the ground and top level parameter estimators,
$(\hat{\delta}_{12},\hat{\delta}_{23})$ and $\hat{\delta}_{13|2}$,
respectively, that is, the ratios between the variances of the ML and
alternative estimators in question. In a Gumbel copula, the dependence
increases with the parameter $\delta$. Kendall's $\tau$ is $0$ when
$\delta= 1$ and tends to $1$ as $\delta\rightarrow\infty$. The
three estimators are rather efficient in general, with
IFM on top, followed by SP and finally SSP. As the true margins
are known, this is not that surprising, and agrees with the results
of Kim \textit{et~al.}~\cite{silvapulle07}. All three estimators lose asymptotic
efficiency with increasing dependence at the ground level,
that is, for $\delta_{12}$ and $\delta_{23}$, whereas SP and SSP gain
efficiency at the top level. The asymptotic variances of all three estimators
actually decrease with increasing dependence at both levels, though
not as fast as for ML. As expected, SSP is overall less efficient than
IFM and SP, but the difference is quite small at
the top level.
\end{Ex}
\begin{Ex}\label{ex:student}
Consider the five-dimensional D-vine with Student's \textit{t}-margins and
Student's \textit{t}-copulae:
\[
f_{12345}(x_{1},x_{2},x_{3},x_{4},x_{5};\B{\nu}_{M},\B{\rho},\B
{\nu}_{C}) =  \prod_{l=1}^{5}f(x_{l};\nu_{l})\prod_{j=1}^{4}\prod
_{i=1}^{5-j}c(u_{i|v_{ij}},u_{i+j|v_{ij}};\rho_{i,i+j|v_{ij}},\nu
_{i,i+j|v_{ij}}),\vadjust{\goodbreak}
\]
with $\B{\nu}_{M}=(\nu_{1},\ldots,\nu_{5})$,
$\B{\rho}=(\rho_{12},\ldots,\rho_{15|234})$,
$\B{\nu}_{C}=(\nu_{12},\ldots,\nu_{15|234})$, and
\begin{eqnarray*}
f(x;\nu) &= & \frac{\Gamma((\nu+1)/2 )}{\sqrt{\uppi\nu
}\Gamma(\nu/2 )}\biggl (1+\frac{x^{2}}{\nu} \biggr)^{-(\nu
+1)/2},\\
 c(u_{i},u_{j};\rho,\nu) &= & \Gamma\biggl(\frac{\nu+2}{2}
\biggr)\Gamma\biggl(\frac{\nu}{2} \biggr)\biggl (1+\frac{t_{\nu}^{-1}(u_{i})^{2}}{\nu}
\biggr)^{(\nu+1)/2} \biggl(1+\frac{t_{\nu}^{-1}(u_{j})^{2}}{\nu}
\biggr)^{(\nu+1)/2}\\
&&{}\Big/\biggl(\Gamma\biggl(\frac{\nu+1}{2} \biggr)^{2}\sqrt{1-\rho
^{2}} \\
&&\hphantom{\Big/\biggl(}{}\times\biggl(1+\frac{t_{\nu}^{-1}(u_{i})^{2}+t_{\nu}^{-1}(u_{j})^{2}-2\rho
t_{\nu}^{-1}(u_{i})t_{\nu}^{-1}(u_{j})}{\nu(1-\rho^{2})} \biggr)^{(\nu+2)/2}\biggr),\\
 u_{i|v_{ij}} &= &
h(u_{i|v_{i,j-1}},u_{i+j-1|v_{i,j-1}};\rho_{i,i+j-1|v_{i,j-1}},\nu
_{i,i+j-1|v_{i,j-1}}),\\
 u_{i+j|v_{ij}} &= &
h(u_{i+j|v_{i+1,j-1}},u_{i+1|v_{i+1,j-1}};\rho
_{i+1,i+j|v_{i+1,j-1}},\nu_{i+1,i+j|v_{i+1,j-1}}),\\
 h(u,v;\rho,\nu)
&= & t_{\nu+1} \bigl(\sqrt{\nu+1} \bigl(t_{\nu}^{-1}(u)-\rho t_{\nu}^{-1}(v)
\bigr) \bigl( \bigl(\nu+t_{\nu}^{-1}(v)^{2} \bigr)(1-\rho^{2}) \bigr)^{-1/2} \bigr),
\end{eqnarray*}
with $u_{i}=t_{\nu}(x_{i})$, $t_{\nu}$ being the c.d.f. of the Student's
\textit{t}-distribution with $\nu$ degrees of freedom. This is a five-dimensional
extension of the example model considered in Sections~\ref{subsubsec:sp}
and~\ref{subsec:ssp}, that is, none of the copulae share parameters,
nor do the margins. The number of parameters is therefore~$25$.

\begin{table}
\tablewidth=\textwidth
\tabcolsep=0pt
  \caption{Asymptotic relative efficiencies of $\hat{\B{\rho}}$ and
$\hat{\B{\nu}}_{C}$ from Example \protect\ref{ex:student}, averaged over
each level, for various parameter sets}\label{tab:expstud}
\begin{tabular*}{\textwidth}{@{\extracolsep{\fill}}lllllllll@{}}
\hline
 & \multicolumn{2}{l}{Level 1} &
\multicolumn{2}{l}{Level 2} & \multicolumn{2}{l}{Level 3} &
\multicolumn{2}{l}{Level 4}\\[-5pt]
& \multicolumn{2}{l}{\hrulefill} &
\multicolumn{2}{l}{\hrulefill} & \multicolumn{2}{l}{\hrulefill} &
\multicolumn{2}{l@{}}{\hrulefill}\\
 & $\hat{\rho}$ & $\hat{\nu}$ & $\hat
{\rho}$ & $\hat{\nu}$ & $\hat{\rho}$ & $\hat{\nu}$ & $\hat{\rho
}$ & $\hat{\nu}$   \\
\hline
\multicolumn{9}{c}{\textsc{IFM}}\\
$(\rho,\nu)=(0.3,6)$ & 0.988 & 0.996 & 0.988 & 0.997 & 0.998 & 0.996
& 0.997 & 0.998\\
 $(\rho,\nu)=(0.7,6)$ & 0.935 & 0.913 & 0.961 & 0.988
& 0.984 & 0.979 & 0.968 & 0.996\\
 $(\rho,\nu)=(0.3,20)$ & 0.997 & 0.996
& 0.993 & 0.999 & 0.990 & 0.998 & 0.997 & 0.999\\
 $(\rho,\nu)=(0.7,20)$
& 0.952 & 0.992 & 0.962 & 0.993 & 0.969 & 0.988 & 0.991 & 0.989 \\
\multicolumn{9}{c}{\textsc{SP}} \\
$(\rho,\nu)=(0.3,6)$ & 0.952 & 0.985 & 0.973 & 0.987 & 0.991 & 0.992
& 0.998 & 0.997\\
 $(\rho,\nu)=(0.7,6)$ & 0.872 & 0.883 & 0.915 & 0.956
& 0.974 & 0.965 & 0.988 & 0.977 \\
$(\rho,\nu)=(0.3,20)$ & 0.965 & 0.963
& 0.994 & 0.970 & 0.989 & 0.991 & 0.983 & 0.992\\
 $(\rho,\nu)=(0.7,20)$
& 0.938 & 0.926 & 0.966 & 0.992 & 0.958 & 0.985 & 0.994 & 0.990\\
\multicolumn{9}{c}{\textsc{SSP}} \\
$(\rho,\nu)=(0.3,6)$ & 0.890 & 0.907 & 0.932 & 0.937 & 0.938 & 0.985
& 0.946 & 0.996\\
 $(\rho,\nu)=(0.7,6)$ & 0.852 & 0.855 & 0.861 & 0.934
& 0.925 & 0.948 & 0.967 & 0.959\\
 $(\rho,\nu)=(0.3,20)$ & 0.941 & 0.955
& 0.992 & 0.964 & 0.981 & 0.975 & 0.973 & 0.974 \\
$(\rho,\nu)=(0.7,20)$
& 0.870 & 0.911 & 0.950 & 0.990 & 0.968 & 0.981 & 0.980 & 0.985 \\
\hline
\end{tabular*}
\end{table}

In this case, it is infeasible to compute the asymptotic
covariance matrices, both analytically and numerically. Therefore,
we resort to simulation and Monte Carlo methods. More specifically,
we have generated $N=250$ samples of size $n=10\mbox{,}000$ from the above
distribution with four different parameter sets. For each sample, we have
estimated the PCC parameters $\B{\rho}$ and $\B{\nu}_{C}$ using the
four estimators. Finally, we have computed the sample covariance
matrices of the resulting estimates. The four parameter sets we have
considered are $(\rho,\nu)=(0.3,6),(0.7,6),(0.3,20),(0.7,20)$, where
we let $\rho_{12}=\cdots=\rho_{15|234}=\rho$, $\nu_{12}=\cdots=\nu
_{15|234}=\nu$,
fixing the marginal parameters at $\nu_{1}=\cdots=\nu_{5}=6$.
Table~\ref{tab:expstud} shows the resulting relative efficiencies
averaged over each level. The three estimators behave rather
similarly, although IFM once more appears to be the most
efficient, SP the second and SSP the last. More specifically,
their efficiency decreases with increasing dependence (either
higher correlation or lower number of degrees of freedom) at
all levels of the structure. Furthermore, they all become more
efficient with increasing level number. In particular, the SSP
estimator gains with respect to its competitors at the higher
levels, just as for the Gumbel vine in Example~\ref{ex:expgumb}.
Note that an increased efficiency is not synonymous with
a lower estimator variance, but only measures the behaviour
relative to the ML estimator. Actually, the variances of
all four estimators increase with the level of the structure,
as one would expect.
\end{Ex}
\begin{Ex}\label{ex:raindata}
Finally, we have fitted a D-vine to a set of daily precipitation values
recorded from 01.01.1990 to 31.12.2006 at five different meteorological
stations in Norway; Vestby, Ski, L{\o}renskog, Nannestad and Hurdal,
shown on the map (Figure~2,~\ref{sup:Example4.4} of Hob{\ae}k
Haff~\cite{ho2012}). These data were provided by the Norwegian
Meteorological Institute. Moreover, this is one of the data sets
studied in Berg and Aas~\cite{bergaas09}, extended with the series from
L{\o}renskog. We have followed their example and modelled only the
positive precipitation, that is, we have discarded all observations for
which at least one of the stations has recorded zero precipitation,
leaving 2013 out of the original 6209. The aim is to remove the
temporal dependence between the observations, in accordance with our
assumptions. Autocorrelation plots of the resulting data set indicate
that this is reasonable.

Since rain showers tend to be very local, we expect the dependence
between measurements from two proximate stations to be stronger than
from stations that are farther apart. As the stations almost lie on a
straight line (see Figure 2,~\ref{sup:Example4.4} of Hob{\ae}k Haff
\cite{ho2012}), a D-vine ordered according to geography is a very
natural model. More specifically, the chosen dependence structure is
the left-hand side of Figure~\ref{fig:CDvine}, with Vestby, Ski,
L{\o}renskog, Nannestad and Hurdal as variables 1, 2, 3, 4 and 5,
respectively. To find adequate copulae for our structure, we computed
the pseudo-observations, shown in Figure~3 in~\ref{sup:Example4.4} of
Hob{\ae}k Haff~\cite{ho2012}. There are strong indications of upper,
but not of lower tail dependence. We therefore chose Gumbel copulae at
the ground level. An inspection of the data transformed with the
estimated $h$-functions from the preceding level (as described in
Section~\ref{subsec:ssp}) indicated that Gaussian copulae would be
reasonable for the three remaining levels. Finally, according to
histograms of the data (shown in Figure~4,~\ref{sup:Example4.4} of
Hob{\ae}k Haff~\cite{ho2012}), the generalised gamma distribution
(Stacy~\cite{generalisedgamma}) with p.d.f.
\[
f(x;\gamma,\beta,p) = \frac{p}{\beta^{\gamma}\Gamma(\gamma
/p )}x^{\gamma-1}\exp\biggl\{- \biggl(\frac{x}{\beta} \biggr)^{p}\biggr \}
\]
appears to be suitable for the margins. This distribution is gamma
for $p=1$ and exponential if in addition $\gamma=1$. Both the ML and
the IFM estimates of $\gamma$ and $p$ were rather different from $1$,
which confirms that the margins are neither exponential nor gamma
distributions. The actual fitted marginal p.d.f.s are shown in the
histograms of the data (Figure 4).

\begin{table}
\tablewidth=\textwidth
\tabcolsep=0pt
\caption{Estimated parameters with 95$\%$ confidence intervals for the
precipitation data set of Example~\protect\ref{ex:raindata}}\label{tab:raindata}
\begin{tabular*}{\textwidth}{@{\extracolsep{\fill}}llllll@{}}
\hline
Lev. & Par. & \textsc{ML} & \textsc{IFM} &
\textsc{SP} & \textsc{SSP} \\
\hline
 1& $\theta_{12}$ & 4.56 & 4.37 & 4.32 & 4.32\\
 & & (4.44, 4.71) & (4.18, 4.56) & (4.14, 4.50) & (4.14, 4.50)\\
 &  $\theta_{23}$ & 3.02 & 2.92 & 2.91 & 2.90\\
 & & (2.91, 3.13) & (2.80, 3.04) & (2.79, 3.03) & (2.79, 3.03)\\
 & $\theta_{34}$ & 2.53 & 2.47 & 2.47 & 2.47 \\
 & & (2.44, 2.62) & (2.37, 2.57) & (2.37, 2.57) & (2.37, 2.56) \\
 & $\theta_{45}$ & 3.59 & 3.48 & 3.44 & 3.44 \\
 & & (3.45, 3.73) & (3.34, 3.62) & (3.30, 3.58) & (3.30, 3.58)\\[3pt]

2 & $\theta_{13|2}$ & $-$0.17 & $-$0.17 & $-$0.17 &
$-$0.17\\
 & & ($-$0.21, $-$0.13) & ($-$0.21, $-$0.13) & ($-$0.21, $-$0.13) & ($-$0.21, $-$0.13)\\
 & $\theta_{24|3}$ & 0.21 & 0.20 & 0.21 & 0.21 \\
 & & (0.15, 0.27) & (0.16, 0.24) & (0.17, 0.25) & (0.17, 0.25)\\
 & $\theta_{35|4}$ & 0.066 & 0.067 & 0.061 & 0.061 \\
 & & (0.022, 0.11) & (0.031, 0.10) & (0.023, 0.099) & (0.024, 0.098)\\[3pt]

3 & $\theta_{14|23}$ & 0.093 & 0.088 & 0.081 & 0.081\\
 & & (0.055, 0.13) & (0.053, 0.12) & (0.044, 0.12) & (0.044, 0.12)\\
& $\theta_{25|34}$ & 0.050 & 0.043 & 0.033 & 0.033\\
 & & (0.009, 0.091) & (0.008, 0.079) & ($-$0.003, 0.070) & ($-$0.003, 0.070) \\[3pt]
4 & $\theta_{15|234}$ & 0.040 & 0.045 & 0.046 & 0.046\\
 & & (0.006, 0.075) & (0.007, 0.083) & (0.012, 0.080) & (0.012, 0.080)\\
  \hline
\end{tabular*}
\end{table}

We have fitted the described model with each of the four estimators,
using the R-routine \mbox{\texttt{optim()}}. The resulting estimates are
shown in Table~\ref{tab:raindata}, along with 95$\%$
confidence intervals. These were computed by
$\hat{\theta} \pm\Phi^{-1}(0.975)\hat{se}$,
for each of the ten parameters, where $\mathit{\hat{se}}$ is an estimate
of the parameter's asymptotic standard deviation. For the ML
estimator, we computed the sample Fisher matrix
\[
\hat{\I}=\frac{1}{n}\sum_{k=1}^{n}\frac{\partial}{\partial(\B
{\alpha},\B{\theta})}\log(f_{1\ldots5}(x_{1k},\ldots,x_{5k};\hat
{\B{\alpha}}{}^{\mathrm{ML}},\hat{\B{\theta}}{}^{\mathrm{ML}}) ),
\]
where $\B{\alpha}=(\gamma_{1},\beta_{1},p_{1},\ldots,\gamma
_{5},\beta_{5},p_{5})$,
the derivative being calculated numerically. The estimates $\hat{se}$
were then simply the square roots of the diagonal
entries of $\hat{\I}^{(\theta)}$. For the three remaining
estimators, we used parametric bootstrap to obtain $\hat{se}$.
More specifically, we generated $B=500$ bootstrap
samples from
$F_{1\ldots5}(x_{1},\ldots,x_{5};\hat{\B{\alpha}}{}^{\mathrm{IFM}},\hat{\B
{\theta}}{}^{\mathrm{IFM}})$, $C_{1\ldots5}(u_{1},\ldots,u_{5};\hat{\B{\theta
}}{}^{\mathrm{SP}})$
and
$C_{1\ldots5}(u_{1},\ldots,u_{5};\hat{\B{\theta}}{}^{\mathrm{SSP}})$,
estimating the parameters
$\hat{\B{\alpha}}{}^{\mathrm{IFM},b},\break\hat{\B{\theta}}{}^{\mathrm{IFM},b},
\hat{\B{\theta}}{}^{\mathrm{SP},b},\hat{\B{\theta}}{}^{\mathrm{SSP},b}$,
$b=1,\ldots,B$. Finally, we let the $\hat{se}$s be the sample
standard deviations of the bootstrap estimates.

At the ground level, the parameter estimates are overall high.
This indicates a strong positive dependence between large
amounts of precipitation in stations that are close in distance,
as anticipated. The IFM, SP and SSP estimators give similar
values. However, the ML estimates are rather different,
though the 95$\%$ confidence intervals overlap with the other
estimators'. As noted earlier, this indicates that the chosen
univariate margins or copulae are not quite adequate. Since the
SP and SSP estimates are virtually the same, the problem
appears to be the margins.

The second level models the conditional dependencies of two
stations that are separated by one, given the one between them.
All four estimators agree that this conditional dependence
is negative between Vestby and L{\o}renskog, positive between
the pair Ski and Nannestad, whereas L{\o}renskog and Hurdal
are almost conditionally independent. At the top two levels,
the estimated copulae are close to the independence copula,
as expected. Actually, the SP and SSP confidence intervals
indicate that the copula $C_{25|34}$ is not significantly
different from independence, which can be an important aspect
for practical purposes.
\end{Ex}
\section{Concluding remarks}\label{sec:conclusion}
There are various estimators for the parameters of a pair-copula
construction, among those the stepwise semiparametric estimator,
which is designed for this particular dependence structure. Although
previously suggested, it has never been formally introduced. In
this paper, we have presented its asymptotic properties, as well as
the estimation algorithm for the two most common types of PCCs,
namely D- and C-vines.

Compared to alternatives such as maximum likelihood,
inference functions for margins and semiparametric estimation,
SSP is in general asymptotically less efficient. The SSP
estimator has a higher variance than the alternatives. Nonetheless,
the loss of efficiency is rather low, and decreases with the
construction level, as shown in a couple of examples. For the
set of five precipitation series, the SSP estimates are actually
almost indistinguishable from the SP ones. Moreover, the SSP
estimator is semiparametrically so for the Gaussian copula. To
compare the alternative estimators' performance more thoroughly,
we plan to perform a large simulation study.

One of the main advantages of the SSP estimator, is that it is
computationally tractable even in high dimensions, as opposed
to its competitors. Moreover, it provides start values required
by the other estimators. Finally, determining the pair-copulae
of a PCC is a stepwise procedure, that involves parameter estimates
from preceding levels. The SSP estimator lends itself perfectly
to that task.

For simplicity, we have only considered C- and D-vines. Equivalent
results are, however, easily obtained for the more general class
of regular vines. Further, we have partitioned the parameter vector
into marginal and dependence parameters. This excludes some
distributions, such as the multivariate Student's \textit{t}. However, if
one does not need the flexibility to model the margins and
dependence structure separately, as well different types of
dependence between the various pairs of variables, a PCC is
unnecessarily complex. Moreover, we have assumed the observations
to be independent, identically distributed. In practice, the
parameter estimation often includes a preliminary step to deal
with deviations from these assumptions (Chen and Fan~\cite{chenfan06}), for
instance GARCH filtration of time series data. The effect of such
an additional step on the SSP estimator is a subject for future
work.
\begin{appendix}\label{sec:appendix}
\section*{Appendix}
\subsection{\texorpdfstring{Matrices $\B{\K}_{\theta}$ and $\B{\J}_{\theta}$}
{Matrices K theta and J theta}}\label{subsec:sspcomp}

As stated in Section~\ref{subsec:ssp}, the matrices
$\B{\K}_{\theta} = \E(\B{\psi}_{\theta}\B{\psi}_{\theta}^{T} )$
and
$\B{\J}_{\theta} = \E(-\B{\psi}_{\theta,\theta} )$
are block diagonal and block lower triangular, respectively, that is,
$\B{\K}_{\theta,i,j}=\B{0}, i\neq j$ and
$\B{\J}_{\theta,i,j}=\B{0},\break i < j$. This follows from the
structure of the $\psi$-functions, as shown below.

We start with $\B{\J}_{\theta,i,j}$, where $i < j$. Then,
$\B{\J}_{\theta,i,j} = \E(-\frac{\partial^{2} \psi
_{i}(u_{1},\ldots,u_{d};\B{\theta}_{1},\ldots,\B{\theta
}_{i})}{\partial\B{\theta}_{i}\,\partial\B{\theta}_ {j}^{T}} ),$
with $ \psi_{i}$ from \eqref{eqn:psi-funcD}.
Since none of the copulae at level $i$ are functions of the
parameters at a following level $j$,
$\frac{\partial\psi_{i}(u_{1},\ldots,u_{d};\B{\theta}_{1},\ldots
,\B{\theta}_{i})}{\partial\B{\theta}_ {j}} = \B{0}$.
Hence, $\B{\J}_{\theta,i,j}=\B{0}, i < j$.

Assume now that $i < j$, and let
$\mathbf{u}=(u_{1},\ldots,u_{d})=(\mathbf{u}_{w_{ki}},\mathbf{u}_{-w_{ki}})$.
Then,
\begin{eqnarray*}
\B{\K}_{\theta,i,j} & = &\E\biggl( \biggl(\frac{\partial\psi_{i}(u_{1},\ldots
,u_{d};\B{\theta}_{1},\ldots,\B{\theta}_{i})}{\partial\B{\theta
}_{i}} \biggr) \biggl(\frac{\partial\psi_{j}(u_{1},\ldots,u_{d};\B{\theta
}_{1},\ldots,\B{\theta}_{j})}{\partial\B{\theta}_ {j}} \biggr)^{T}\biggr )\\
 &=& \int_{\mathbf{u}}\frac{\partial}{\partial\B{\theta}_{i}}\sum
_{k=1}^{d-i}\log c_{k,k+i|v_{ki}}\frac{\partial}{\partial\B{\theta
}_{j}^{T}}\sum_{l=1}^{d-j}\log c_{l,l+j|v_{lj}}c_{1\ldots d}\,\mathrm{d}\mathbf{u}\\
 &=& \sum_{k=1}^{d-i}\sum_{l=1}^{d-j}\int_{\mathbf{u}_{w_{ki}}}\frac
{\partial}{\partial\B{\theta}_{i}}\log c_{k,k+i|v_{ki}}\int_{\mathbf
{u}_{-w_{ki}}}\frac{1}{c_{l,l+j|v_{lj}}}\frac{\partial}{\partial\B
{\theta}_{j}^{T}}c_{l,l+j|v_{lj}}c_{1\ldots d}\,\mathrm{d}\mathbf{u}_{-w_{ki}}\,\mathrm{d}\mathbf
{u}_{w_{ki}}.
\end{eqnarray*}
Under the conditions of Theorem~\ref{theo:ssp}, we may exchange the
integration and differentiation in the inner integral. Thus,
\begin{eqnarray*}
\B{\K}_{\theta,i,j} & =& \sum_{k=1}^{d-i}\sum_{l=1}^{d-j}\int_{\mathbf
{u}_{w_{ki}}}\frac{\partial}{\partial\B{\theta}_{i}}\log
c_{k,k+i|v_{ki}}\frac{\partial}{\partial\B{\theta}_{j}^{T}} \biggl(\int
_{\mathbf{u}_{-w_{ki}}}\frac{c_{l,l+j|v_{lj}}c_{1\ldots
d}}{c_{l,l+j|v_{lj}}}\,\mathrm{d}\mathbf{u}_{-w_{ki}} \biggr)\,\mathrm{d}\mathbf{u}_{w_{ki}}\\
 & =& \sum
_{k=1}^{d-i}\sum_{l=1}^{d-j}\int_{\mathbf{u}_{w_{ki}}}\frac{\partial
}{\partial\B{\theta}_{i}}\log c_{k,k+i|v_{ki}}\frac{\partial
}{\partial\B{\theta}_{j}^{T}} \biggl(\int_{\mathbf{u}_{-w_{ki}}}c_{1\ldots
d}\,\mathrm{d}\mathbf{u}_{-w_{ki}} \biggr)\,\mathrm{d}\mathbf{u}_{w_{ki}}\\
 & = &\sum_{k=1}^{d-i}\sum
_{l=1}^{d-j}\int_{\mathbf{u}_{w_{ki}}}\frac{\partial}{\partial\B
{\theta}_{i}}\log c_{k,k+i|v_{ki}}\frac{\partial}{\partial\B
{\theta}_{j}^{T}}c_{w_{ki}}\,\mathrm{d}\mathbf{u}_{w_{ki}}.
\end{eqnarray*}
The pair-copulae composing $c_{w_{ki}}$, situated in levels
$1,\ldots,i$, are not functions of parameters from a following
level $j$. Thus, $\frac{\partial}{\partial\B{\theta
}_{j}}c_{w_{ki}}=\B{0}$.
Consequently, $\B{\K}_{\theta,i,j}=\B{0}, i < j$. The exact same
argument can be repeated for $i > j$. Hence,
$\B{\K}_{\theta,i,j}=\B{0}, i \neq j$.\vspace*{-3pt}

\subsection{\texorpdfstring{Proof of Theorem \protect\ref{theo:sspnorm}}
{Proof of Theorem 2}}\label{subsec:proofsspnorm}\vspace*{-6pt}
\begin{pf}
In two dimensions, the SSP estimator is the same as the
SP estimator, which was shown to be semiparametrically efficient
by Klaassen and Wellner~\cite{klaassen97}. In three dimensions, we have computed
the asymptotic covariance matrices for comparison. As shown
in Example~\ref{ex:normdist}, the covariance matrices of the
SP and SSP estimators, $\mathbf{V}^{\mathrm{SP}}$ and $\mathbf{V}^{\mathrm{SSP}}$,
respectively, are equal. Thus, the SSP estimator is
semiparametrically efficient also for the three-dimensional
Gaussian copula.

Assume now that it is true for the $(d-1)$-dimensional Gaussian
copula. Further, for the $d$-dimensional model, partition the
covariance matrix $\mathbf{V}^{\mathrm{SP}}$ as
$\mathbf{V}_{1}^{\mathrm{SP}}=\mathbf{V}_{1\ldots d-2,1\ldots d-2}^{\mathrm{SSP}}$,
$\mathbf{V}_{12}^{\mathrm{SP}}=\mathbf{V}_{1\ldots d-2,d-1}^{\mathrm{SSP}}$ and
$V_{2}^{\mathrm{SP}}=V_{d-1,d-1}^{\mathrm{SSP}}$, and likewise for $\mathbf{V}^{\mathrm{SSP}}$,
$\mathbf{V}^{\mathrm{ML}}$, $\mathbf{B}^{\mathrm{SP}}$, $\mathbf{B}^{\mathrm{SSP}}$, $\B{\I}_{\theta}$,
$\B{\I}^{(\theta)}$ and $\B{\J}_{\theta}$.
As the SP estimator is semiparametrically efficient, $\mathbf{V}^{\mathrm{ML}}$
for the Gaussian copula must be the same, regardless of the margins.
Moreover, when all margins are normal, $\hat{\B{\theta}}{}^{\mathrm{ML}}$ is
simply the empirical correlation matrix. Adding an extra dimension
leaves the remaining estimators unchanged. Hence, $\mathbf{V}_{1}^{\mathrm{ML}}$,
corresponding to the $(d-1)$-dimensional sub-model, will be the same
as for the $(d-1)$-dimensional Gaussian copula. The same argument
can repeated for all $(d-1)$-dimensional sub-models, covering all
levels but the top. Due to its levelwise structure, the SSP estimator
for a given sub-model is unaffected when adding an extra dimension,
and so must the corresponding block of $\mathbf{V}^{\mathrm{SSP}}$ be. Accordingly,
we must have $\mathbf{V}_{1}^{\mathrm{SSP}}=\mathbf{V}_{1}^{\mathrm{SP}}=\mathbf{V}_{1}^{\mathrm{ML}}$.
Hence, it remains to show that $\mathbf{V}_{12}^{\mathrm{SSP}}=\mathbf{V}_{12}^{\mathrm{SP}}$
and $V_{2}^{\mathrm{SSP}}=V_{2}^{\mathrm{SP}}$, related to the estimators
$\hat{\theta}{}^{\mathrm{SP}}_{1d|v_{1d}}$ and $\hat{\theta}{}^{\mathrm{SSP}}_{1d|v_{1d}}$
for the top level copula. According to Theorem 1 from Tsukahara~\cite{tsukahara05}
and Theorem~\ref{theo:ssp}, respectively,
\[
\sqrt{n} (\hat{\theta}{}^{\mathrm{SP}}_{1d|v_{1d}}-\theta_{1d|v_{1d}} )
\stackrel{d}{\longrightarrow}
Z_{\mathrm{SP}} \sim\N(0,V_{2}^{\mathrm{SP}} )
\]
and
\[
\sqrt{n} (\hat{\theta}{}^{\mathrm{SSP}}_{1d|v_{1d}}-\theta_{1d|v_{1d}} )
\stackrel{d}{\longrightarrow}
Z_{\mathrm{SSP}} \sim\N(0,V_{2}^{\mathrm{SSP}} ),
\]
as $n \rightarrow\infty$. Now, define
$\mathbf{U}_{n}= (\mathbf{U}_{n1},\ldots,\mathbf{U}_{nn} )$,
with $\mathbf{U}_{nj}=(F_{1n}(X_{1j}),\ldots,F_{dn}(X_{dj})),\break j=1,\ldots,d$,
and let
\begin{eqnarray*}
\B{\Psi}^{\mathrm{SP}}(\mathbf{U}_{n};\hat{\B{\theta}}{}^{\mathrm{SP}}) & =& \frac
{1}{n}\sum_{k=1}^{n}\B{\phi}^{\mathrm{SP}} (\mathbf{U}_{n};\hat{\B{\theta
}}{}^{\mathrm{SP}} )=\B{0},\\[-2pt]
 \B{\Psi}^{\mathrm{SSP}}(\mathbf{U}_{n};\hat{\B{\theta}}{}^{\mathrm{SSP}})
& = &\frac{1}{n}\sum_{k=1}^{n}\B{\phi}^{\mathrm{SSP}} (\mathbf{U}_{n};\hat{\B
{\theta}}{}^{\mathrm{SSP}} )=\B{0},\vadjust{\goodbreak}
\end{eqnarray*}
be the estimating equations of the SP and SSP estimators,
respectively. Further, let
\[
\Psi(\mathbf{U}_{n};\B{\theta})  = \Psi_{d(d-1)/2}^{\mathrm{SP}}(\mathbf
{U}_{n};\B{\theta}) = \Psi_{d(d-1)/2}^{\mathrm{SSP}}(\mathbf{U}_{n};\B
{\theta}) = \frac{1}{n}\sum_{k=1}^{n}\frac{\partial}{\partial
\theta_{1d|v_{1d}}}\log(c_{1\ldots d}(\mathbf{U}_{n};\B{\theta}) ).
\]
According to Theorem 1 from Tsukahara~\cite{tsukahara05},
\begin{eqnarray*}
\Psi(\mathbf{U}_{n};\hat{\B{\theta}}{}^{\mathrm{SP}} ) &= & \Psi(\mathbf{U}_{n};\B
{\theta})+\frac{\partial\Psi(\mathbf{U}_{n};\B{\theta})}{\partial
\theta_{1d|v_{1d}}} (\hat{\theta}{}^{\mathrm{SP}}_{1d|v_{1d}}-\theta
_{1d|v_{1d}} )\\
&&{} +  \frac{\partial\Psi(\mathbf{U}_{n};\B{\theta
})}{\partial\B{\theta}_{1\rightarrow d-2}^{T}} (\hat{\B{\theta
}}{}^{\mathrm{SP}}_{1\rightarrow d-2}-\B{\theta}_{1\rightarrow d-2} )+\mathrm{o}_{P}
\biggl(\frac{1}{n}\biggr ) = 0.
\end{eqnarray*}
Likewise, using Theorem~\ref{theo:ssp}, one obtains
\begin{eqnarray*}
\Psi(\mathbf{U}_{n};\hat{\B{\theta}}{}^{\mathrm{SSP}} ) &= & \Psi(\mathbf{U}_{n};\B
{\theta})+\frac{\partial\Psi(\mathbf{U}_{n};\B{\theta})}{\partial
\theta_{1d|v_{1d}}} (\hat{\theta}{}^{\mathrm{SSP}}_{1d|v_{1d}}-\theta
_{1d|v_{1d}} )\\
&&{} +  \frac{\partial\Psi(\mathbf{U}_{n};\B{\theta
})}{\partial\B{\theta}_{1\rightarrow d-2}^{T}} (\hat{\B{\theta
}}{}^{\mathrm{SSP}}_{1\rightarrow d-2}-\B{\theta}_{1\rightarrow d-2} )+\mathrm{o}_{P}
\biggl(\frac{1}{n}\biggr ) = 0.
\end{eqnarray*}
Hence,
\[
\sqrt{n} (\hat{\theta}{}^{\mathrm{SSP}}_{1d|v_{1d}}-\hat{\theta
}{}^{\mathrm{SP}}_{1d|v_{1d}} ) = \frac{\mathbf{A}_{1}}{A_{2}}\sqrt{n} (\hat{\B
{\theta}}{}^{\mathrm{SP}}_{1\rightarrow d-2}-\hat{\B{\theta}}{}^{\mathrm{SSP}}_{1\rightarrow
d-2} )+\mathrm{o}_{P}
\biggl(\frac{1}{n}\biggr ),
\]
with
\[
\mathbf{A}_{1}=\frac{\partial\Psi(\mathbf{U}_{n};\B{\theta})}{\partial\B
{\theta}_{1\rightarrow d-2}^{T}}=\frac{1}{n}\sum_{k=1}^{n}\frac
{\partial^{2}}{\partial\theta_{1d|v_{1d}} \,\partial\B{\theta
}_{1\rightarrow d-2}^{T}}\log(c_{1\ldots d}(\mathbf{U}_{n};\B{\theta}) )
\]
and
\[
A_{2}=\frac{\partial\Psi(\mathbf{U}_{n};\B{\theta})}{\partial\theta
_{1d|v_{1d}}}=\frac{1}{n}\sum_{k=1}^{n}\frac{\partial^{2}}{\partial
\theta_{1d|v_{1d}}^{2}} \log(c_{1\ldots d}(\mathbf{U}_{n};\B{\theta}) ).
\]
According to the assumption,
\[
\sqrt{n} (\hat{\B{\theta}}{}^{\mathrm{SSP}}_{1\rightarrow d-2}-\B{\theta
}_{1\rightarrow d-2} ) \stackrel{d}{\longrightarrow} \mathbf{Y} \sim\N
_{d(d-1)/2-1} (\B{0},\mathbf{V}_{1}^{\mathrm{ML}} ),\qquad n \rightarrow\infty.
\]
Thus,
\begin{eqnarray*}
\sqrt{n} (\hat{\B{\theta}}{}^{\mathrm{SP}}_{1\rightarrow d-2}-
\hat{\B{\theta}}{}^{\mathrm{SSP}}_{1\rightarrow d-2} )
& =& \sqrt{n} (\hat{\B{\theta}}{}^{\mathrm{SP}}_{1\rightarrow d-2}-
\B{\theta}_{1\rightarrow d-2} )-\sqrt{n}(\hat{\B{\theta}}{}^{\mathrm{SSP}}_{1\rightarrow d-2}-
\B{\theta}_{1\rightarrow d-2} )\\
& \stackrel{p}{\longrightarrow}& \B{0},\qquad n
\rightarrow\infty.
\end{eqnarray*}
Moreover, under the assumed conditions,
$\mathbf{A}_{1} \stackrel{p}{\longrightarrow} -\B{\I}_{\theta,12}^{T}$
and
$A_{2} \stackrel{p}{\longrightarrow} -\I_{\theta,2}$,
as $n \rightarrow\infty$. Hence,
\[
\sqrt{n} (\hat{\theta}{}^{\mathrm{SSP}}_{1d|v_{1d}}-\theta_{1d|v_{1d}}^{\mathrm{SP}} )
\stackrel{p}{\longrightarrow} 0,\vadjust{\goodbreak}
\]
which means that $Z_{\mathrm{SP}} \stackrel{d}{=} Z_{\mathrm{SSP}}$. In other words,
$V_{2}^{\mathrm{SSP}}=V_{2}^{\mathrm{SP}}$.
Moreover,
\begin{eqnarray*}
\mathbf{V}_{12}^{\mathrm{SSP}}& = & \frac{1}{\I_{\theta,2}} (-\mathbf{V}_{1}^{\mathrm{ML}}\B
{\I}_{\theta,12}+\B{\J}_{\theta,12}^{-1}\mathbf{B}_{12}^{\mathrm{SSP}} ),\\
V_{2}^{\mathrm{SSP}} &= & \frac{1}{\I_{\theta,2}} \biggl(1+\frac{B_{2}^{\mathrm{SSP}}}{\I
_{\theta,2}}+ \frac{1}{\I_{\theta,2}}\B{\I}_{\theta,12}^{T}\mathbf
{V}_{1}^{\mathrm{ML}}\B{\I}_{\theta,12}-\frac{2}{\I_{\theta,2}}\B{\I
}_{\theta,12}^{T}\B{\J}_{\theta,12}^{-1}\mathbf{B}_{12}^{\mathrm{SSP}}\biggr ) \\
&= &
\frac{1}{\I_{\theta,2}^{2}} (\I_{\theta,2}+B_{2}^{\mathrm{SSP}}-\B{\I
}_{\theta,12}^{T}\mathbf{V}_{1}^{\mathrm{ML}}\B{\I}_{\theta,12} )-2\B{\I
}_{\theta,12}^{T}\mathbf{V}_{12}^{\mathrm{SSP}}.
\end{eqnarray*}
Correspondingly for SP,
\begin{eqnarray*}
\mathbf{V}_{12}^{\mathrm{SP}}& = & \frac{1}{\I_{\theta,2}} \biggl(-\mathbf{V}_{1}^{\mathrm{ML}}\B{\I
}_{\theta,12}+\B{\I}_{1}^{(\theta)} \biggl(\mathbf{B}_{12}^{\mathrm{SP}}-\frac
{B_{2}^{\mathrm{SP}}}{\I_{\theta,2}}\B{\I}_{\theta,12} \biggr)\biggr ), \\
V_{2}^{\mathrm{SP}} &= &
\frac{1}{\I_{\theta,2}} \biggl(1+\frac{B_{2}^{\mathrm{SP}}}{\I_{\theta,2}}+\frac
{1}{\I_{\theta,2}}\B{\I}_{\theta,12}^{T}\mathbf{V}_{1}^{\mathrm{ML}}\B{\I
}_{\theta,12}-\frac{2}{\I_{\theta,2}}\B{\I}_{\theta,12}^{T}\B
{\I}_{1}^{(\theta)} \biggl(\mathbf{B}_{12}^{\mathrm{SP}}-\frac{B_{2}^{\mathrm{SP}}}{\I_{\theta
,2}}\B{\I}_{\theta,12} \biggr) \biggr)\\
 &= & \frac{1}{\I_{\theta,2}^{2}} (\I
_{\theta,2}+B_{2}^{\mathrm{SP}}-\B{\I}_{\theta,12}^{T}\mathbf{V}_{1}^{\mathrm{ML}}\B{\I
}_{\theta,12} )-2\B{\I}_{\theta,12}^{T}\mathbf{V}_{12}^{\mathrm{SP}}.
\end{eqnarray*}
Since the estimating equation for $\theta_{1d|v_{1d}}$ is
the same for SP and SSP, $B_{2}^{\mathrm{SSP}}=B_{2}^{\mathrm{SP}}$. Moreover,
$V_{2}^{\mathrm{SSP}}=V_{2}^{\mathrm{SP}}$. Consequently, $\mathbf{V}_{12}^{\mathrm{SSP}}=\mathbf{V}_{12}^{\mathrm{SP}}$.
\end{pf}
\subsection{\texorpdfstring{Covariance matrices from Example \protect\ref{ex:normdist}}
{Covariance matrices from Example 4.1}}\label{subsec:normexcomp}
The asymptotic covariance matrix of the ML estimator
is given by
\[
\mathbf{V}^{\mathrm{ML}}  = \frac{1}{2}
\pmatrix{
2(1-\rho_{12}^{2})^{2} & v_{13} & v_{23}\cr
 v_{13} & 2(1-\rho
_{23}^{2})^{2} & v_{12}\cr
 v_{23} & v_{12} & 2(1-\rho_{13}^{2})^{2}
},
\]
with $v_{ik} = 2\rho_{ik}(1-\rho_{il}^{2})(1-\rho_{lk}^2)-\rho
_{il}\rho_{lk}|\mathbf{R}|$.
For the SSP estimator, we have
\[
\mathbf{V}^{\mathrm{SSP}}  = \B{\J}_{\theta}^{-1}\B{\K}_{\theta} (\B{\J
}_{\theta}^{-1} )^{T}+\B{\J}_{\theta}^{-1}\mathbf{B}_{\theta}^{\mathrm{SSP}}
(\B{\J}_{\theta}^{-1} )^{T},
\]
where
{\fontsize{9.5pt}{\baselineskip}\selectfont{
\begin{eqnarray*}
\B{\K}_{\theta} & =&
\pmatrix{
\displaystyle\frac{1+\rho_{12}^{2}}{(1-\rho_{12}^{2})^{2}} & \displaystyle\frac
{k_{12}}{(1-\rho_{12}^{2})(1-\rho_{23}^{2})} & 0\cr
 \displaystyle\frac{k_{12}}{(1-\rho_{12}^{2})(1-\rho_{23}^{2})} & \displaystyle\frac{1+\rho
_{23}^{2}}{(1-\rho_{23}^{2})^{2}} & 0 \cr
0 & 0 & \displaystyle\frac{|\mathbf{R}|+2(\rho
_{13}-\rho_{12}\rho_{23})^2}{|\mathbf{R}|^{2}}
},\\
 \B{\J}_{\theta} & =&
\pmatrix{
\displaystyle\frac{1+\rho_{12}^{2}}{(1-\rho_{12}^{2})^{2}} & 0 & 0 \cr
0 & \displaystyle\frac
{1+\rho_{23}^{2}}{(1-\rho_{23}^{2})^{2}} & 0 \cr
\displaystyle\frac{j_{23}}{|\mathbf
{R}|^{2}} & \displaystyle\frac{j_{12}}{|\mathbf{R}|^{2}} & \displaystyle\frac{|\mathbf{R}|+2(\rho
_{13}-\rho_{12}\rho_{23})^2}{|\mathbf{R}|^{2}}
},\\
 \mathbf{B}_{\theta}^{\mathrm{SSP}} & =&
\pmatrix{
\displaystyle\frac{\rho_{12}^{2}(1+\rho_{12}^{2})}{(1-\rho_{12}^{2})^{2}} &\displaystyle\frac{\rho_{23}b_{12}+\rho_{12}b_{23}-\rho_{12}\rho_{23}a}{2(1-\rho_{12}^{2})(1-\rho_{23}^{2})}& \displaystyle\frac{(\rho_{13}-\rho_{12}\rho_{23})b_{12}}{2(1-\rho_{12}^{2})|\mathbf{R}|}\cr
\displaystyle\frac{\rho_{23}b_{12}+\rho_{12}b_{23}-\rho_{12}\rho_{23}a}{2(1-\rho_{12}^{2})(1-\rho_{23}^{2})} & \displaystyle\frac{\rho_{23}^{2}(1+\rho_{23}^{2})}{(1-\rho_{23}^{2})^{2}}&\displaystyle\frac{(\rho_{13}-\rho_{12}\rho_{23})b_{23}}{2(1-\rho_{23}^{2})|\mathbf{R}|}\cr
\displaystyle\frac{(\rho_{13}-\rho_{12}\rho_{23})b_{12}}{2(1-\rho_{12}^{2})|\mathbf{R}|} & \displaystyle\frac{(\rho_{13}-\rho_{12}\rho_{23})b_{23}}{2(1-\rho_{23}^{2})|\mathbf{R}|}&\displaystyle\frac{(1+\rho_{13}^{2})(\rho_{13}-\rho_{12}\rho_{23})^2}{|\mathbf{R}|^{2}}},
\end{eqnarray*}}}%
with
$a = 1+\rho_{12}^{2}+\rho_{13}^{2}+\rho_{23}^{2}$,
$k_{12} = (\rho_{13}-\rho_{12}\rho_{23})(|\mathbf{R}|+\rho
_{13}^{2}-\rho_{12}^{2}\rho_{23}^{2})$,
$j_{ik} = -\rho_{ik}|\mathbf{R}|+2(\rho_{il}-\rho_{lk}\rho_{ik})(\rho
_{lk}-\rho_{il}\rho_{ik})$
and
$b_{ik} = \rho_{ik}a (1-\frac{2}{1-\rho_{ik}^{2}} )+2(1+\rho
_{ik}^{2})\frac{\rho_{ik}+\rho_{il}\rho_{lk}}{1-\rho_{ik}^{2}}$.
The resulting covariance matrix is
\[
\mathbf{V}^{\mathrm{SSP}}  = \frac{1}{2}
\pmatrix{
2(1-\rho_{12}^{2})^{2} & v_{13} & v_{23} \cr
v_{13} & 2(1-\rho
_{23}^{2})^{2} & v_{12}\cr
 v_{23} & v_{12} & 2(1-\rho_{13}^{2})^{2}}
= \mathbf{V}^{\mathrm{ML}}.\vspace*{-6pt}
\]
\end{appendix}

\section*{Acknowledgements}\label{sec:acknowledgements}
This work is funded by Statistics for Innovation, (sfi)$^2$.
I thank my supervisors Arnoldo Frigessi and Kjersti Aas for
very helpful discussions and comments. I also thank the referees
and Associate Editor for their help to improve this paper
with their good comments and suggestions. Finally, I would like
to give special thanks to Hideatsu Tsukahara for having clarified
the validity of Theorem 1 in Tsukahara~\cite{tsukahara05} for
general $m$.\vspace*{-6pt}

%
\begin{supplement}
\sname{Supplement A}\label{sup:SSPalg}
\stitle{SSP estimation algorithms for D- and C-vines\\}
\slink[doi]{10.3150/12-BEJ413SUPPA}
\sdatatype{.pdf}
\sfilename{BEJ413\_suppa.pdf}
\sdescription{Estimation algorithms for the stepwise
semiparametric estimator for D- and C-vines.\vspace*{-3pt}}
\end{supplement}
\begin{supplement}
\sname{Supplement B}\label{sup:Example4.4}
\stitle{Figures and table from Example~\ref{ex:raindata}}
\slink[doi]{10.3150/12-BEJ413SUPPB}
\sdatatype{.pdf}
\sfilename{BEJ413\_suppb.pdf}
\sdescription{Figures 2, 3 and 4, as well as Table 4
from Example~\ref{ex:raindata}.\vspace*{-6pt}}
\end{supplement}
%

\printhistory

\end{document}